\documentclass[10pt]{article} 
\usepackage[letterpaper, bindingoffset=0.2in, left=.4in, right=1in, top=1in, bottom=1in, footskip=.25in]{geometry}
\usepackage{blindtext}
\usepackage{multicol}
\usepackage{lipsum}
\usepackage{cuted}
\usepackage{mwe}
\usepackage{amsmath, amsthm, amssymb} % May not all be necessary
\usepackage{graphicx} % For including pictures
\usepackage{algorithm}
\usepackage{algorithmic}
\usepackage{wrapfig}
\usepackage{listings}
\usepackage[colorlinks=true, urlcolor=blue, citecolor=black, linkcolor=black]{hyperref} % For formatting (clickable) URLs
\usepackage{subcaption} % For captioning multi-panel figures
\makeatletter
\def\th@definition{\thm@notefont{}\normalfont}
\makeatother
\theoremstyle{definition}
\newtheorem{fact}{Fact}
\newtheorem{theorem}{Theorem}
\newcommand{\cupdot}{\mathbin{\mathaccent\cdot\cup}}

\title{The Smooth Power of the ``Neandertal Method''}
\author{
{Aaron Montag\textsuperscript{1}, Tim Reinhardt\textsuperscript{2}, Jürgen Richter-Gebert\textsuperscript{3}}\\
\vspace{2pt}
\textsuperscript{}\quad\\
Department of Mathematics, CIT, Technical University of Munich \\
\textsuperscript{1} {\tt montag@posteo.net}\quad
\textsuperscript{2} {\tt tim.reinhardt@tum.de}\quad
\textsuperscript{3} {\tt richter@tum.de}\\
}
\date{}	% Suppress any date on submissions

\begin{document}

\maketitle

% Prevent page number 1 from being printed on the first page.
\thispagestyle{empty}

\begin{abstract}
We describe an algorithmic method to transform a Euclidean wallpaper pattern into a \emph{Circle Limit}-style picture à la Escher. The design goals for the method are to be mathematically sound, aesthetically pleasing and fast to compute. It turns out that a certain class of conformal maps is particularly well-suited for the problem. Moreover, in our specific application, a very simple method – sometimes jokingly called the ``Neandertal method'' for its almost brutal simplicity – proves to be highly efficient, as it can easily be parallelized to be run on the GPU, unlike many other approaches.
\end{abstract}
\bigskip
\begin{multicols}{2}
%\raggedright
% Bridges papers are usually no more than 8 pages in length. So
% there's really no need to have numbered sections, unless the
% author really needs to refer to sections by number within the paper's text.
% So to suppress sequential section numbers, append an asterisk to 
% the \section command, as in:

%%%%%%%%%%%%%%%%%%%%%%%%%%%%%%%%%%%%%%%%%%%
\section{\texorpdfstring{\raggedright Hyperbolization of Euclidean Tilings}{Hyperbolization of Euclidean Tilings}}
``I worked terribly hard to finally finish that litho, and then with gritted teeth, spent another four days making beautiful prints of that extremely complex circle limit in colors. Each print is a series of twenty printings: five pieces, and each piece four times. All this with the remarkable feeling that this work is a milestone in my development, and that nobody, except myself, will ever realize this.'', M.\ C.\ Escher wrote in 1960 in a letter to his son about producing his \emph{Circle Limit} prints.

%\begin{wrapfigure}{r}{0.45\textwidth}
%\vspace{-.7cm}
%\begin{center}
%\includegraphics[width=0.43\textwidth]{Images/Escher3.jpg}
%\end{center}
%\caption{A close-up of Eschers's ``Circle Limit III''}
%\end{wrapfigure}

%\begin{figure*}[t]
% \centering
% \begin{subfigure}[b]{0.46\textwidth}
% 	\centering
%	\includegraphics[width=.95\textwidth]{Images/Escher3.jpg}
%	\caption{Close-up view of ``Circle Limit III''}
%	\label{EscherCL3}
% \end{subfigure}
% \hfill
% \begin{subfigure}[b]{0.46\textwidth}
% 	\centering
%	\includegraphics[width=.95\textwidth]{Images/BirdFischTurtle.png}
% 	\caption{``Fish, Duck, Turtle'' (Symmetry Drawing No.\ 69)}
% 	\label{EscherFDT}
% \end{subfigure}
%\caption{Two works by M.\ C.\ Escher}
%\label{EscherOverview}
%\end{figure*}

\begin{figure}[H]
 \centering
 %\begin{subfigure}[b]{0.46\textwidth}
 	\centering
	\includegraphics[width=.45\textwidth]{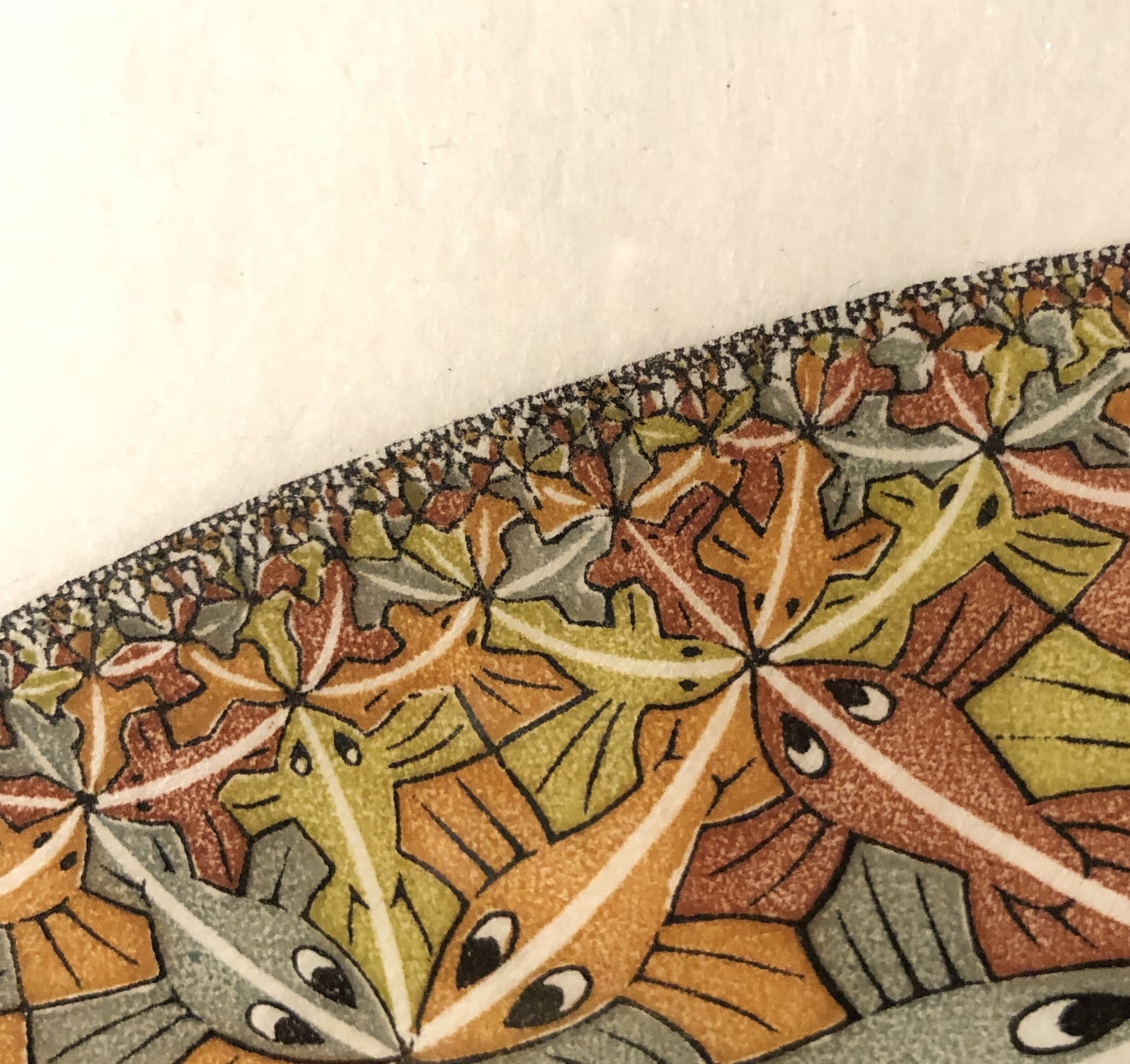}
	\caption{Close-up view of ``Circle Limit III''}
\label{EscherCL3}
 %\end{subfigure}
\end{figure}

\begin{figure}[H]
 \centering
 %\begin{subfigure}[b]{0.46\textwidth}
 	\centering
	\includegraphics[width=.45\textwidth]{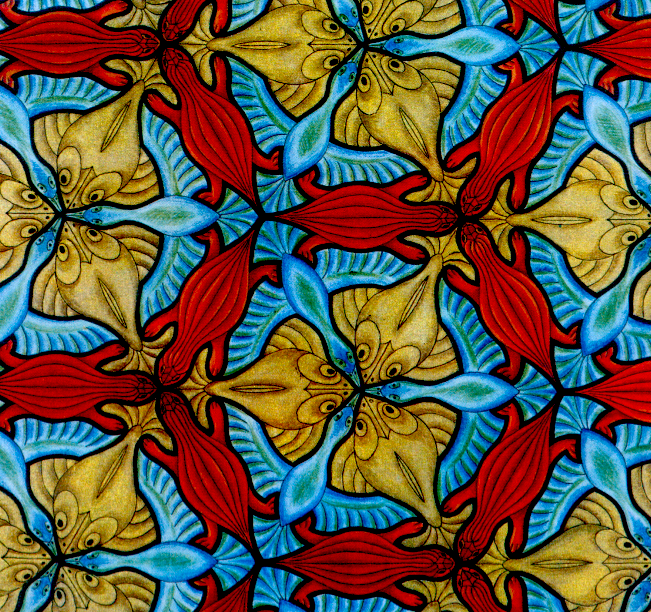}
	\caption{``Fish, Duck, Turtle'' (Symmetry Drawing No. 69)}
	\label{EscherFDT}
 %\end{subfigure}
\end{figure}

A close look at any of Escher's \emph{Circle Limits} (e.g., Figure~\ref{EscherCL3}) immediately reveals his mastery of fine detail. These works mark some of his most structurally and artistically advanced efforts to capture the concept of \emph{infinity} within the finite bounds of a print.
The \emph{Circle Limit} series is based on regular tilings in the hyperbolic plane, an idea Escher encountered through brief contact with H.~S.~M.~Coxeter.
His ability to grasp the essence of a graphical concept without requiring full formal background allowed him to create artworks that arguably did more to popularize and communicate basic ideas of hyperbolic geometry than most mathematics books on that topic.
In contrast to the small number – only four – of \emph{Circle Limit} prints, the quantity of his Euclidean tessellations (over 100, e.g., Figure~\ref{EscherFDT}) is overwhelming. Birds, boats, crabs, fish, horses, and lizards are just some of the motifs featured in Escher's studies of Euclidean tilings.

The project described in this article originated over fifteen years ago, when the third author and Martin von Gagern worked on algorithmic methods to transform Euclidean ornamental patterns into hyperbolic ones; see~\cite{HYP}.
The idea emerged from a fortunate coincidence: We were simultaneously developing a hyperbolic ornament doodling program and software that could automatically detect the symmetry type of an ornament from a scan and extract its fundamental cell.
Combining both led to an automatic \emph{hyperbolization} program capable of transforming many of Escher's Euclidean tessellations into \emph{Circle Limit}-style patterns.
Compared to previous work by Jos Leys~\cite{Leys} and Douglas Dunham~\cite{Dunham} who also attempted to create hyperbolizations of Escher tilings, one of our main goals was to find the mathematically most satisfying method of hyperbolization.
Surprisingly, some approaches even provide a quasi-unique way to create a hyperbolic version of a Euclidean ornament and are therefore preferable over others.
We have since advanced those principles into a parallelized GPU algorithm.
In this article, we report on our current method, which we believe is mathematically satisfying, easy to implement, aesthetically pleasing, and capable of producing smooth real-time transitions between different symmetry types. Moreover, the algorithm has the educational potential to communicate the concept of conformal deformations by providing an extremely simple method that can create such deformations in specific situations.

Besides their direct applicability to ornamental patterns, such algorithms for conformal deformations can – and have – also been used in several purely mathematical contexts, like the description of configuration spaces of linkages or parameterizations of the Clebsch cubic.

\section{\texorpdfstring{\raggedright Ornaments in the Euclidean and in the Hyperbolic Plane}{Ornaments in the Euclidean and in the Hyperbolic Plane}}
An \emph{ornament} is an image in either the Euclidean or hyperbolic plane that exhibits symmetries. Formally, we call a non-constant function $f \colon \mathbb{R}^2 \to C$ a \emph{Euclidean ornament}, where $C$ is a suitable colorspace such as $C=[0,1]$ for grayscale or $C=[0,1]^3$ for RGB color.
A \emph{symmetry} of $f$ is an isometry $\phi \colon \mathbb{R}^2 \to \mathbb{R}^2$ such that $f(x) = f(\phi(x))$ holds for all $x\in\mathbb{R}^2$. The set of all symmetries of $f$ forms a group, the \emph{symmetry group} of the ornament.
\emph{Hyperbolic ornaments} are defined analogously, replacing $\mathbb{R}^2$ by a model of hyperbolic geometry and using hyperbolic isometries.

\begin{figure*}[t]
	\centering
	\setlength{\fboxsep}{0pt}
	\setlength{\fboxrule}{0.5pt}
	\begin{subfigure}[b]{0.09\textwidth}
		\captionsetup{labelformat=empty}
		\fbox{\includegraphics[width=0.95\textwidth]{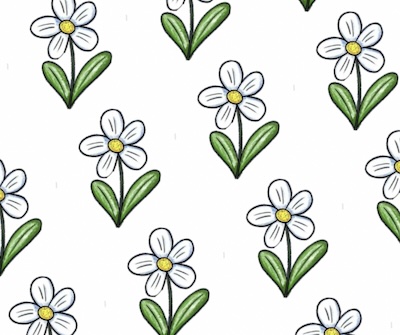}}
		\caption{\centering \makebox[1.0\linewidth][c]{$\mathbf{p1} \mid \mathbf{\bigcirc}$}}
	\end{subfigure}
	\hfill
	\begin{subfigure}[b]{0.09\textwidth}
		\captionsetup{labelformat=empty}
		\fbox{\includegraphics[width=0.95\textwidth]{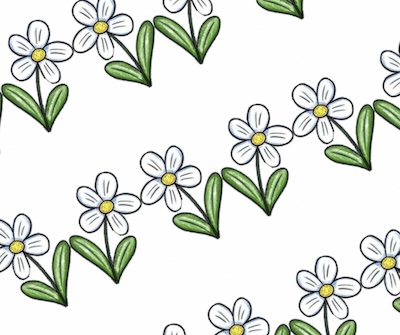}}
		\caption{\centering \makebox[1.0\linewidth][c]{$\mathbf{pm} \mid \mathbf{{\ast}{\ast}}$}}
	\end{subfigure}
	\hfill
	\begin{subfigure}[b]{0.09\textwidth}
		\captionsetup{labelformat=empty}
		\fbox{\includegraphics[width=0.95\textwidth]{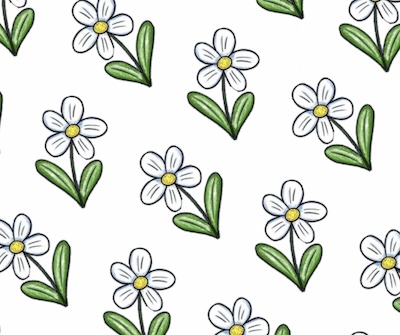}}
		\caption{\centering \makebox[1.0\linewidth][c]{$\mathbf{pg} \mid \mathbf{{\times}{\times}}$}}
	\end{subfigure}
	\hfill
	\begin{subfigure}[b]{0.09\textwidth}
		\captionsetup{labelformat=empty}
		\fbox{\includegraphics[width=0.95\textwidth]{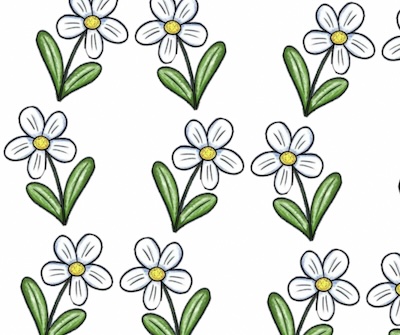}}
		\caption{\centering \makebox[1.0\linewidth][c]{$\mathbf{cm} \mid \mathbf{{\ast}{\times}}$}}
	\end{subfigure}
	\hfill
	\begin{subfigure}[b]{0.09\textwidth}
		\captionsetup{labelformat=empty}
		\fbox{\includegraphics[width=0.95\textwidth]{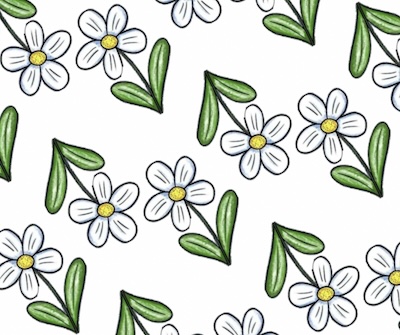}}
		\caption{\centering \makebox[1.0\linewidth][c]{$\mathbf{p2} \mid \mathbf{2222}$}}
	\end{subfigure}
	\hfill
	\begin{subfigure}[b]{0.09\textwidth}
		\captionsetup{labelformat=empty}
		\fbox{\includegraphics[width=0.95\textwidth]{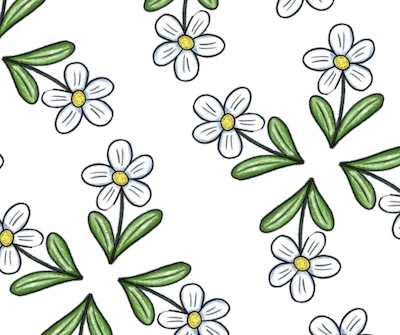}}
		\caption{\centering \makebox[1.0\linewidth][c]{$\mathbf{pmm} \mid \mathbf{{\ast}2222}$}}
	\end{subfigure}
	\hfill
	\begin{subfigure}[b]{0.09\textwidth}
		\captionsetup{labelformat=empty}
		\fbox{\includegraphics[width=0.95\textwidth]{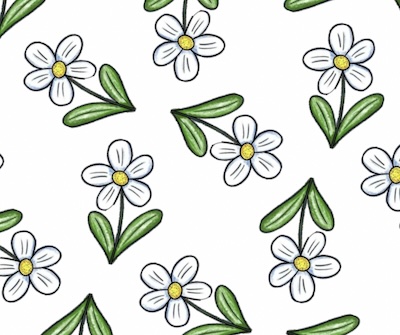}}
		\caption{\centering \makebox[1.0\linewidth][c]{$\mathbf{pgg} \mid \mathbf{22{\times}}$}}
	\end{subfigure}
	\hfill
	\begin{subfigure}[b]{0.09\textwidth}
		\captionsetup{labelformat=empty}
		\fbox{\includegraphics[width=0.95\textwidth]{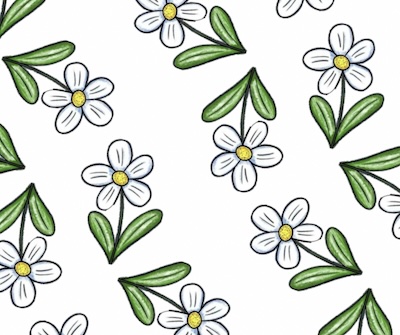}}
		\caption{\centering \makebox[1.0\linewidth][c]{$\mathbf{pmg} \mid \mathbf{22{\ast}}$}}
	\end{subfigure}
	\hfill
	\begin{subfigure}[b]{0.09\textwidth}
		\captionsetup{labelformat=empty}
		\fbox{\includegraphics[width=0.95\textwidth]{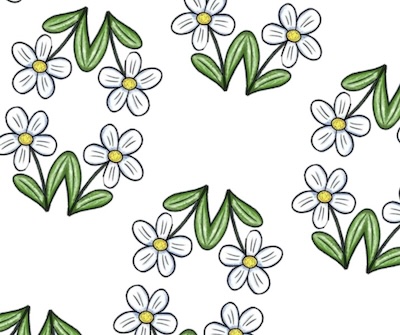}}
		\caption{\centering \makebox[1.0\linewidth][c]{$\mathbf{cmm}\; {\mid}\; \mathbf{2{\ast}22}$}}
	\end{subfigure}
	\\[5pt]
	\begin{subfigure}[b]{0.09\textwidth}
		\captionsetup{labelformat=empty}
		\fbox{\includegraphics[width=0.95\textwidth]{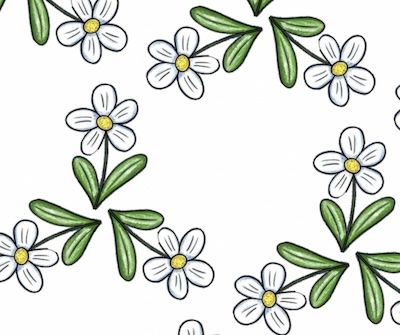}}
		\caption{\centering \makebox[1.0\linewidth][c]{$\mathbf{p3} \mid \mathbf{333}$}}
	\end{subfigure}
	\hfill
	\begin{subfigure}[b]{0.09\textwidth}
		\captionsetup{labelformat=empty}
		\fbox{\includegraphics[width=0.95\textwidth]{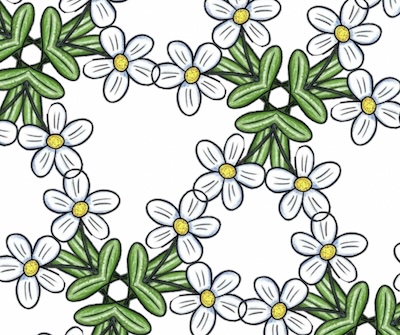}}
		\caption{\centering \makebox[1.0\linewidth][c]{$\mathbf{p3m1}\; {\mid} \;\mathbf{{\ast}333}$}}
	\end{subfigure}
	\hfill
	\begin{subfigure}[b]{0.09\textwidth}
		\captionsetup{labelformat=empty}
		\fbox{\includegraphics[width=0.95\textwidth]{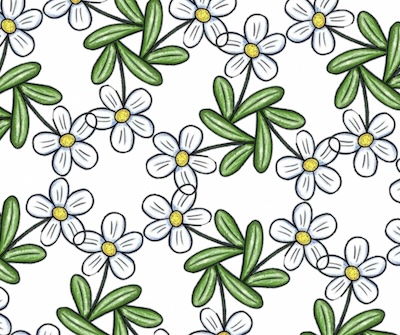}}
		\caption{\centering \makebox[1.0\linewidth][c]{$\mathbf{p31m}\; {\mid}\; \mathbf{3{\ast}3}$}}
	\end{subfigure}
	\hfill
	\begin{subfigure}[b]{0.09\textwidth}
		\captionsetup{labelformat=empty}
		\fbox{\includegraphics[width=0.95\textwidth]{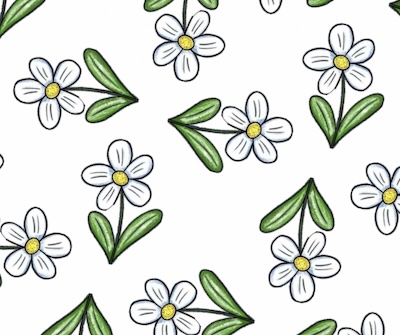}}
		\caption{\centering \makebox[1.0\linewidth][c]{$\mathbf{p4} \mid \mathbf{442}$}}
	\end{subfigure}
	\hfill
	\begin{subfigure}[b]{0.09\textwidth}
		\captionsetup{labelformat=empty}
		\fbox{\includegraphics[width=0.95\textwidth]{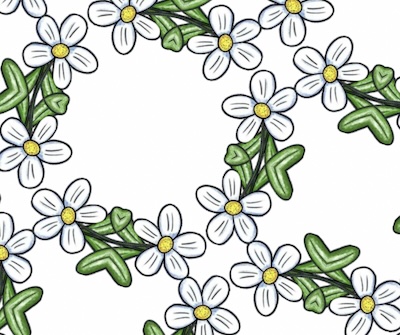}}
		\caption{\centering \makebox[1.0\linewidth][c]{$\mathbf{p4m}\, {\mid}\, \mathbf{{\ast}442}$}}
	\end{subfigure}
	\hfill
	\begin{subfigure}[b]{0.09\textwidth}
		\captionsetup{labelformat=empty}
		\fbox{\includegraphics[width=0.95\textwidth]{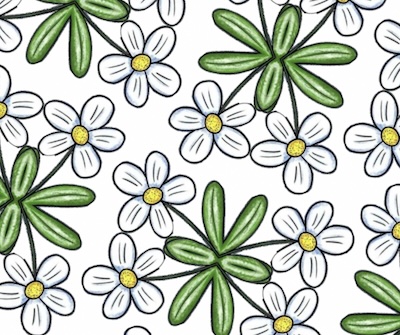}}
		\caption{\centering \makebox[1.0\linewidth][c]{$\mathbf{p4g} \mid \mathbf{4{\ast}2}$}}
	\end{subfigure}
	\hfill
	\begin{subfigure}[b]{0.09\textwidth}
		\captionsetup{labelformat=empty}
		\fbox{\includegraphics[width=0.95\textwidth]{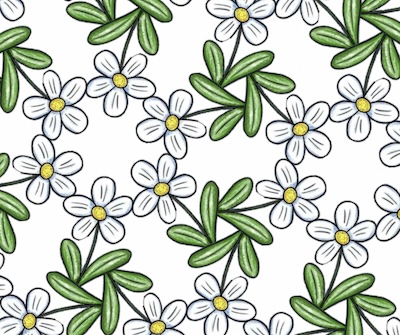}}
		\caption{\centering \makebox[1.0\linewidth][c]{$\mathbf{p6} \mid \mathbf{632}$}}
	\end{subfigure}
	\hfill
	\begin{subfigure}[b]{0.09\textwidth}
		\captionsetup{labelformat=empty}
		\fbox{\includegraphics[width=0.95\textwidth]{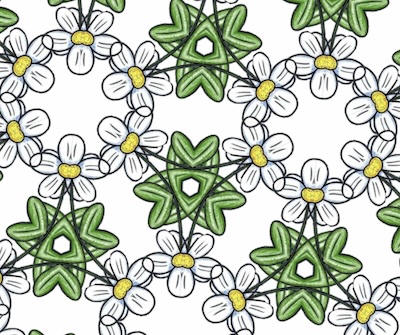}}
		\caption{\centering \makebox[1.0\linewidth][c]{$\mathbf{p6m}\,{\mid}\, \mathbf{{\ast}632}$}}
	\end{subfigure}
	\caption{The 17 wallpaper groups}
	\label{gr17}
\end{figure*}

In the Euclidean plane, there are four types of isometries: translations, rotations, reflections, and glide reflections.
A \emph{wallpaper pattern} is a Euclidean ornament with translational symmetries in two linearly independent directions. As a consequence, every wallpaper ornament is determined by a compact region without non-trivial symmetries – the \emph{fundamental cell} – whose repeated copies tile the entire plane.
In this case, a classical theorem of geometric group theory states that there are exactly 17 non-trivial symmetry groups: the \emph{crystallographic} or \emph{wallpaper} groups. They are depicted in Figure~\ref{gr17}, together with two common naming conventions encoding the presence and absence of certain symmetries. The first system comes from classical crystallography and is widely used, the second is mathematically more pleasing.
It is called \emph{orbifold notation} and was introduced by William Thurston (popularized by John~H.~Conway). Roughly speaking, each individual number $n$ indicates the presence of a transitivity class of an $n$-fold rotation center. The symbols $\ast$~and~$\times$ denote the presence of reflectional and glide reflectional symmetries, respectively.
In composed symbols, integers $n$ appearing before $\ast$ represent ``pure'' rotation centers, while those appearing after $\ast$ correspond to rotation centers that lie on a reflection axis.

Our deformation method directly applies to 15 of the 17 groups. For the groups $\mathbf{\bigcirc}$ and $\mathbf{2222}$, it does not apply in general, but it can still be used as long as the Euclidean fundamental region is a rectangle.
We will address the more intricate cases of groups having fewer symmetries in Section~\ref{moreover}.
For simplicity, we first restrict ourselves to groups of types $\mathbf{ijk}$, $\mathbf{{\ast}ijk}$, and $\mathbf{i{\ast}j}$ (the second row of Figure~\ref{gr17}). All of these are based on triangular tilings of the Euclidean plane.
We begin with groups of type $\mathbf{{\ast}ijk}$. After that, we explain how the methods immediately extend to groups of types $\mathbf{ijk}$ and $\mathbf{i{\ast}j}$.

Groups of type $\mathbf{{\ast}ijk}$ have a triangle as the fundamental cell and are generated by reflections only; they correspond to \emph{kaleidoscopes}. To produce a closed ornament of non-overlapping triangles, the angle at each triangle corner must be of the form $180^\circ/n$ for an integer $n$. Moreover, in the Euclidean plane, the angle sum of a triangle equals $180^\circ$.
These constraints allow exactly three triplets of corner angles as solutions: $(60^\circ,60^\circ,60^\circ)$, $(45^\circ,45^\circ,90^\circ)$, and $(30^\circ,60^\circ,90^\circ)$. The integers $n$ are denoted in the orbifold symbols of the corresponding symmetry groups: $\mathbf{{\ast}333}$, $\mathbf{{\ast}442}$, and $\mathbf{{\ast}632}$, respectively.
If the angles were changed so that their sum becomes smaller than $180^\circ$, the fundamental cell can no longer be realized in Euclidean geometry. However, it can be realized in hyperbolic geometry and, as long as all angles are still integer divisors of $180^\circ$, it uniquely extends by reflection over the entire hyperbolic plane without inconsistencies. In that way, a hyperbolic reflection group arises yielding a hyperbolic ornament.
In this article, we use the Poincaré disk model of hyperbolic geometry with respect to a circle $c$. 
In this case, lines correspond to circular arcs which are orthogonal to $c$ and reflections across a line correspond to circle inversions.

\begin{figure*}[t]
	\centering
	\raisebox{-0.5\height}{\includegraphics[width=.33\textwidth]{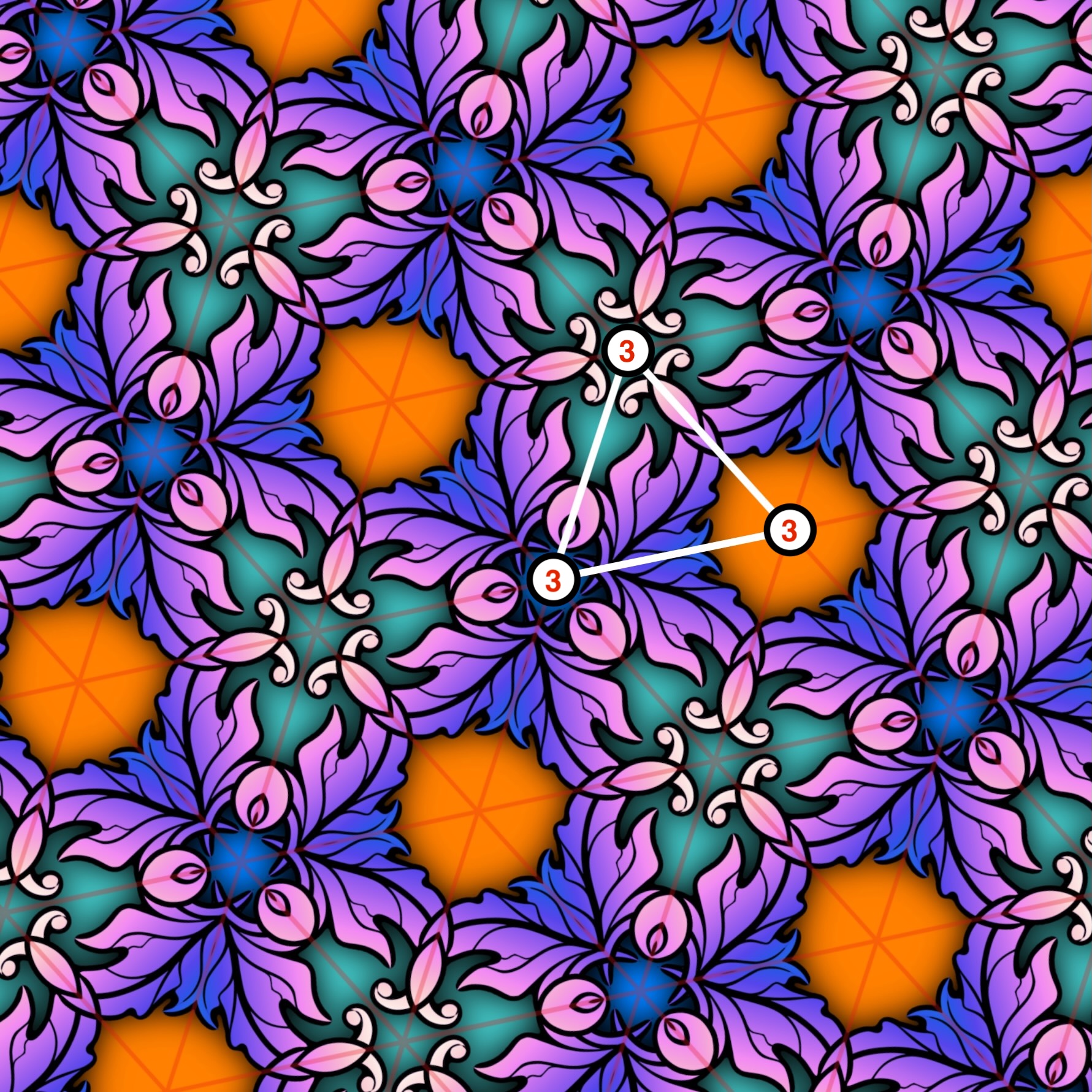}}
	\hspace{4pt}
	\raisebox{-0.5\height}{\huge${\rightarrow}$}
	\hspace{4pt}
	\raisebox{-0.5\height}{\includegraphics[width=.45\textwidth]{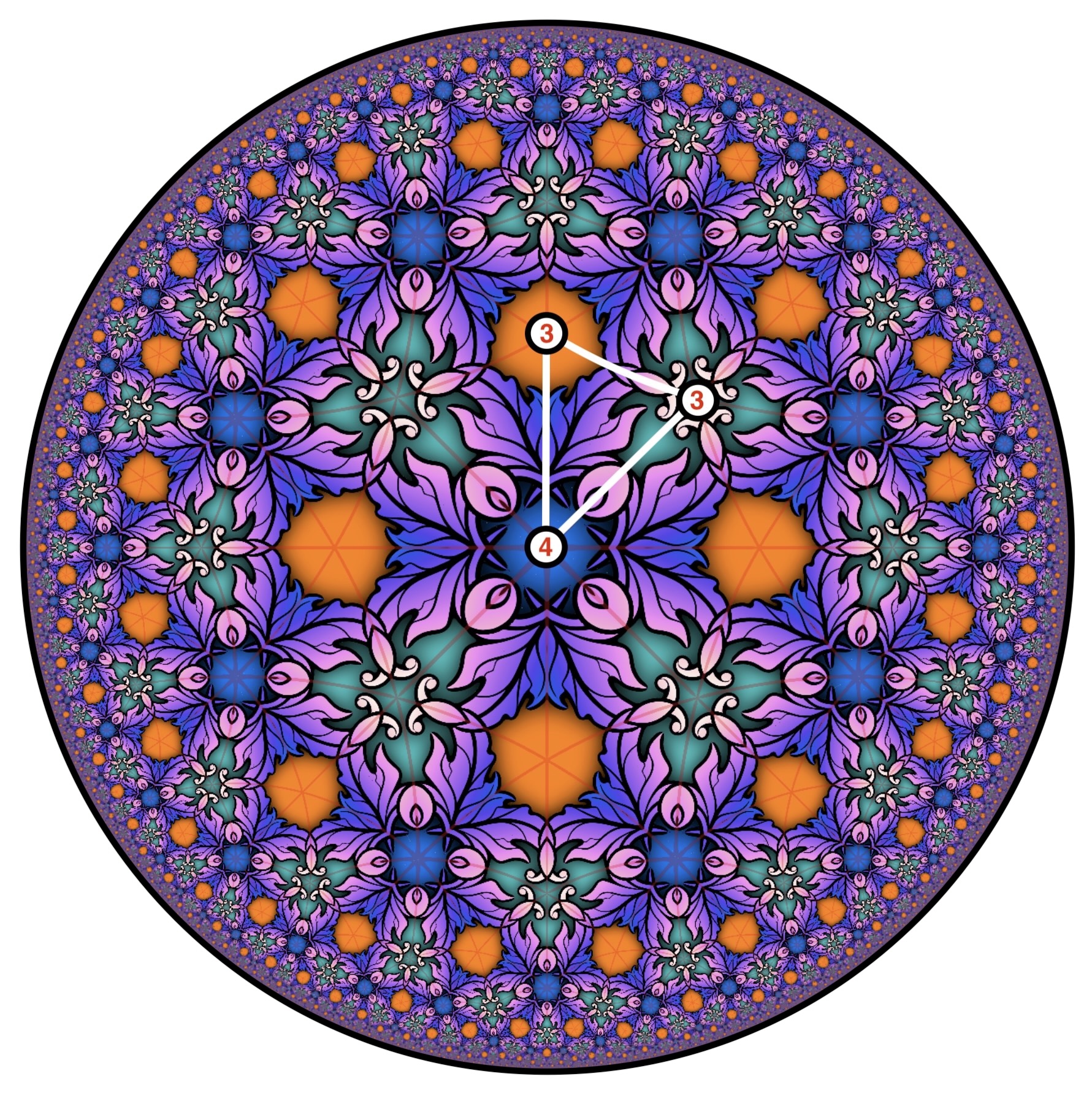}}
	\caption{Hyperbolizing an ornament of type $\mathbf{{\ast}333}$ to type $\mathbf{{\ast}433}$}
	\label{FirstHyperbolization}
\end{figure*}

The principle of hyperbolization is demonstrated in Figure~\ref{FirstHyperbolization}. The left part shows an ornament with symmetry group $\mathbf{{\ast}333}$ whose fundamental cell is marked in white. The entire artistic content is determined by what happens inside this triangle $T_{333}$.
We now transform the pattern into a hyperbolic one by changing the symmetry group to $\mathbf{{\ast}433}$, with the fundamental cell being the hyperbolic triangle $T_{433}$. Ideally, the result should be locally almost indistinguishable from the original picture, while globally creating a new structure. In other words, the artistic content within the fundamental triangle should be disturbed as little as possible. To achieve this, $T_{333}$ has to be mapped bijectively to $T_{433}$ with as minimal distortion as possible.

\section{\texorpdfstring{\raggedright Conformal Maps and the Schwarz Reflection Principle}{Conformal Maps and the Schwarz Reflection Principle}}
\label{sec:ConfAndSRP}
At that point, conformal maps come to the rescue. A map $\mathbb{C}\to \mathbb{C}$ is called \emph{conformal} if its local behaviour is a similarity, i.e.\ a combination of rotation, scaling and translation. Therefore, conformal maps locally preserve angles and ratios of distances, or in other words, they preserve infinitesimal circles and local shapes of figures in general. This makes them extremely suitable in our setting as they locally do not distort the artistic content of an ornament. 
They turn out to be exactly the complex differentiable functions, ensuring that they can be expressed as convergent power series in some neighborhood around each point of their domain.
A~particularly important class of conformal maps are \emph{Möbius transformations}, functions of the form
\[z\mapsto \frac{az+b}{cz+d}\]
for $a,b,c,d \in \mathbb{C}$ satisfying $ad-bc \neq 0$.
Let $D = \{z\in \mathbb{C} \mid |z|<1 \}$ denote the unit disk, considered as a subset of $\mathbb{C}$. For our approach, we will need classical results from complex analysis, which we list below.
 
\begin{fact}[Conformal automorphisms of the unit disk]
\label{confAutom}
Every conformal automorphism of the unit disk $D$ can be expressed as a Möbius transformation. It also maps the unit circle to itself.
\end{fact}

\begin{fact}[Hyperbolic transformations]
\label{hypTrafo}
Identifying the unit disk $D$ with the Poincaré disk model, the conformal automorphisms of $D$ become exactly the orientation-preserving hyperbolic isometries.
\end{fact}
 
\begin{fact}[Three degrees of freedom]
\label{3degrees}
Conformal automorphisms of $D$ have three degrees of freedom. This exactly suffices to move any three points on the unit circle to any three other ones in the same circular order.
\end{fact}
 
\begin{fact}[Riemann mapping theorem]
\label{RMT}
Every region $R$ that is bounded by a Jordan curve can be conformally mapped to the unit disk $D$. This mapping is unique up to a Möbius transformation.
\end{fact}

These four facts are extremely powerful in the context of hyperbolization of ornaments. Recall that our aim for the transformation in Figure~\ref{FirstHyperbolization} is to map a Euclidean triangle $T_{333}$ to a hyperbolic triangle $T_{433}$.
Both triangles are considered as \emph{open} and bounded by a Jordan curve hence by Fact~\ref{RMT} there is a conformal map from $T_{333}$ to the unit disk~$D$. Furthermore, by Fact~\ref{3degrees}, we can find a \emph{unique} mapping $\psi_1$ that sends the corners of $T_{333}$ to a given fixed triangle (say, an equilateral one). Similarly, there exists a second unique conformal map $\psi_2$ that maps $T_{433}$ to the unit disk using the same corner configuration on the unit circle.
The composition
\[\psi \colon T_{433} \to T_{333}, \quad z \mapsto \psi_1^{-1}(\psi_2(z))\]
is a \emph{unique} (!!!) conformal map between the fundamental cells, for a given orientation-preserving assignment of corners.
This solves the problem entirely for the interior of $T_{433}$. We prefer the map $\psi$ with domain $T_{433}$ over its also conformal inverse with domain $T_{333}$ for the following reason.
Recall that the colors of our original ornament are given by $f\colon \mathbb{C} \cong \mathbb{R}^2 \to C$. To assign a color to a specific point $z$ in the fundamental cell $T_{433}$, we can simply map $z \mapsto f(\psi(z))$.

A global conformal hyperbolization of an ornament that extends over the entire plane requires two reasonable ways of extending the map in the fundamental cell to a larger region to coincide. On the one hand, we can extend the map by the rules of symmetry of the two ornaments involved. For a unified treatment, recall that lines can be considered as circles with infinite radius and that circle inversion generalizes ordinary reflections across lines. In $\mathbf{\ast ijk}$, the fundamental cell is bounded by edges that support circular reflection axes in the corresponding symmetry group. Let $\mathcal{I}_c$ denote such a reflection supported by an edge of $T_{433}$, and let $\mathcal{I}_{c_\psi}$ be the corresponding reflection for $T_{333}$. Asking for consistency with respect to symmetry means requiring the operations ``apply the map'' and ``reflect'' to commute.
More formally, we demand:
\begin{equation}
\label{eq:Symmetry}
\psi(\mathcal{I}_c(z))=\mathcal{I}_{c_\psi}(\psi(z))
\end{equation}
for each point $z\in T_{433}$ and each corresponding pair of reflections $(\mathcal{I}_c,\mathcal{I}_{c_\psi})$. We can iteratively apply this consistency criterion to cover the entire pictures in both worlds.

The second concept of extension comes from the fact that conformal maps are analytic.
Hence, the map $\psi$ can be analytically continued from the interior of $T_{433}$ over its boundary as long as no singularity, i.e.\ a rotation center, is hit. Remarkably – and fortunately – these two concepts coincide in our setting.
%This local assignment extends to the global ornament in the most natural way possible.
%Conformal maps are rigid objects: As they are locally given by power series, they can be analytically continued as long as no singularity is hit.
There is a powerful result that makes analytic continuation easy in specific situations: the \emph{Schwarz reflection principle} (SRP). For a circle $c$, we denote its inversion operation by $\mathcal{I}_c$. The following conformal maps extend naturally to the boundaries.
\begin{figure*}[t]
	\centering
	\includegraphics[width=.75\textwidth]{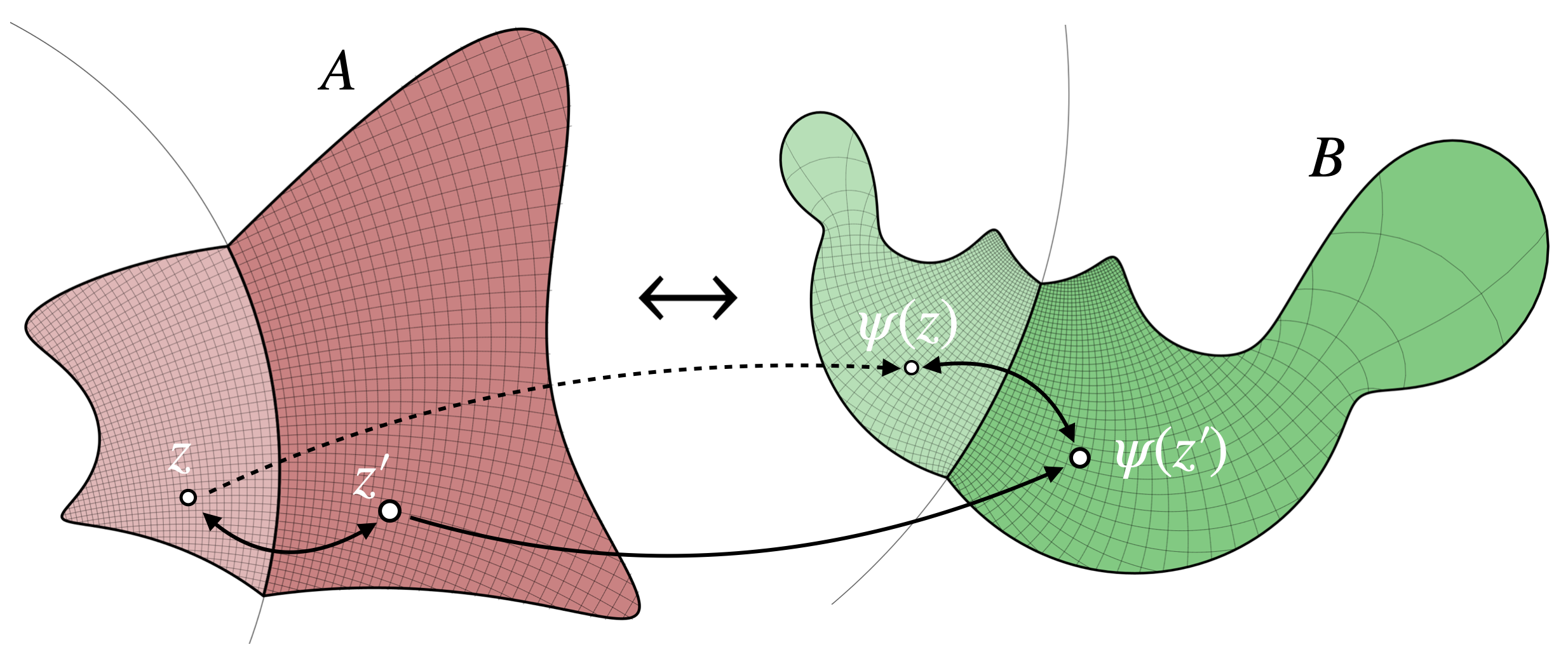}
	\caption{The Schwarz reflection principle}
	\label{SRP}
\end{figure*}

\begin{fact}[Schwarz reflection principle]
\label{factSRP}
Let $A \subset \mathbb{C}$ be a bounded region such that a part $b$ of its boundary is supported by circle $c$. Assume that $\psi\colon A \to B$ is a bijective conformal map such that $\psi(b)$ is supported by a circle $c_\psi$ as well. Then the map $\psi$ can be conformally extended across the boundary $b$. It maps the region $\mathcal{I}_c(A)$ to $\mathcal{I}_{c_\psi}(B)$ according to 
\begin{equation}
\label{eq:SRP}
\psi(z) = \mathcal{I}_{c_\psi}(\psi(\mathcal{I}_c(z))).
\end{equation}
\end{fact}

The situation of Fact~\ref{factSRP} is illustrated in Figure~\ref{SRP}. Suppose that the dark red region $A$ is conformally mapped to the dark green region $B$ such that a part of $\partial A$ is supported by a circle $c$ and gets mapped to a part of $\partial B$ supported by another circle $c_\psi$.
Then the unique analytic extension of $\psi$ maps any point in the light red region $\mathcal{I}_c(A)$ to the light green region by first reflecting it in $c$, applying the original conformal map $\psi$, and finally reflecting in $c_\psi$.

Equation~\ref{eq:SRP} is the same formula as equation~\ref{eq:Symmetry} since reflections are involutions.
Hence, in our situation, the two concepts \emph{extension by symmetry} and \emph{extension by continuation} coincide!

%
%Returning to the example from Figure~\ref{FirstHyperbolization}, we already constructed a conformal map from $T_{433}$ to $T_{333}$. Since the sides of both triangles are circular arcs, we can extend the conformal map beyond its original domain.
%The SRP also tells us exactly how to do this: by a sequence of (circle-)reflections at the sides of the triangles. 

In our example from Figure~\ref{FirstHyperbolization}, we first constructed a conformal map from $T_{433}$ to $T_{333}$. Since the sides of both triangles are circular arcs, we were able to extend the conformal map beyond its original domain by iterated circular reflection.
The SRP ensures that this process yields a map that is globally conformal except for the rotation centers.

\begin{figure*}[t]
	\centering
	\raisebox{-0.5\height}{\includegraphics[width=.37\textwidth]{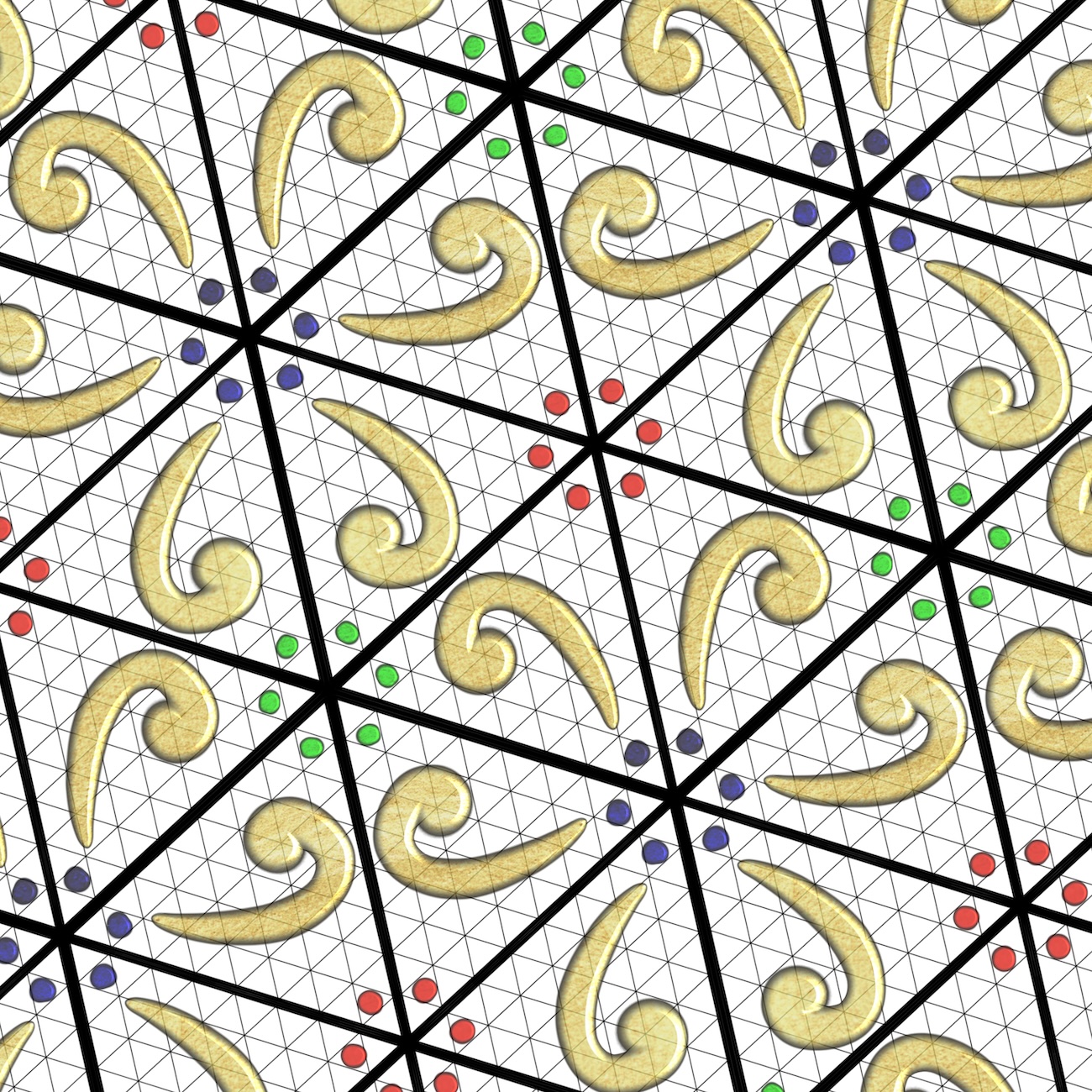}}
	\hspace{4pt}
	\raisebox{-0.5\height}{\huge${\rightarrow}$}
	\hspace{4pt}
	\raisebox{-0.5\height}{\includegraphics[width=.45\textwidth]{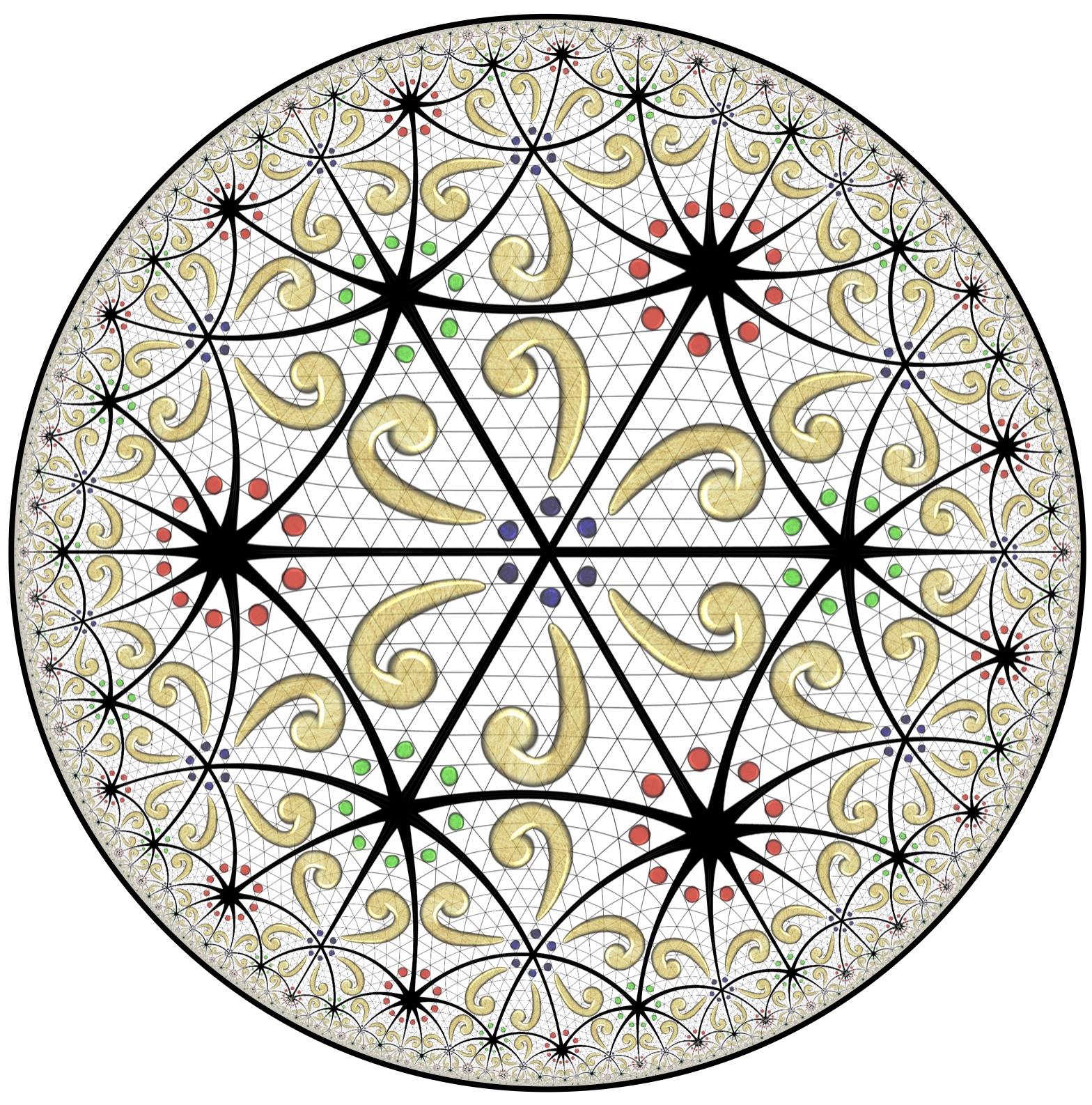}}
	\caption{Hyperbolizing an ornament of type $\mathbf{{\ast}333}$ to type $\mathbf{{\ast}543}$ with corresponding grid lines}
	\label{grid}
\end{figure*}

There is a little subtlety, however. Since the degrees of the rotation centers get altered, a round trip of reflections in the original pattern may not correspond to a round trip of reflections in the image ornament. We will not go into details here but it is exactly the fact that the ornament admits the symmetry together with the fact that the used angles are divisors of $180^\circ$ that ensures that everything is well defined.
Figure~\ref{grid} illustrates the process in a slightly more elaborate example where we transform a symmetry group of type $\mathbf{{\ast}333}$ to type $\mathbf{{\ast}543}$. The hyperbolized fundamental cell is a hyperbolic triangle $T_{543}$ with corner angles 
$(36^\circ,45^\circ,60^\circ)$.
An underlying triangular mesh is shown, illustrating how each triangle is specifically mapped.
Observe how the coordinate lines extend across triangle boundaries to form perfectly smooth, slightly curved lines.

In summary, once the combinatorics is fixed, our method yields a unique conformal map between the fundamental cells. The SRP guarantees that this map \emph{uniquely} extends to an ornament in the entire hyperbolic plane. One may still assume some choice in the position of our triangle $T_{433}$ in the hyperbolic plane. However, choosing another triangle with identical corner angles only results in a hyperbolic motion of the entire picture, as a consequence of Fact~\ref{hypTrafo}.

\section{The Neandertal Method} 
One problem remains: How to compute the conformal map for the fundamental region? Such maps between prescribed shapes are notoriously difficult to construct and require comparably deep insights into complex analysis. Special functions like hypergeometric, Jacobi elliptic or theta functions frequently appear.
In many situations, no explicit analytic expression for the map exists. There is extensive literature on numerical methods for computing conformal mappings – such as the Zipper algorithm and Schwarz-Christoffel transformations – as well as on approaches from discrete differential geometry, including discrete conformal maps and circle packings.
Remarkably, in our situation an extremely simple algorithm produces excellent results. Due to its almost brutal simplicity, it was jokingly called the ``Neandertal method'' in a series of online lectures by Keenan~Crane~\cite{Crane}.
At first glance, the method has two major drawbacks: it requires solving a system of equations with about one million variables and offers little to no flexibility in adapting boundary conditions. Fortunately, both disadvantages can be turned into advantages in our specific application.

\begin{figure}[H]
%\begin{figure}[h!tbp]
	\centering
	\includegraphics[width=.45\textwidth]{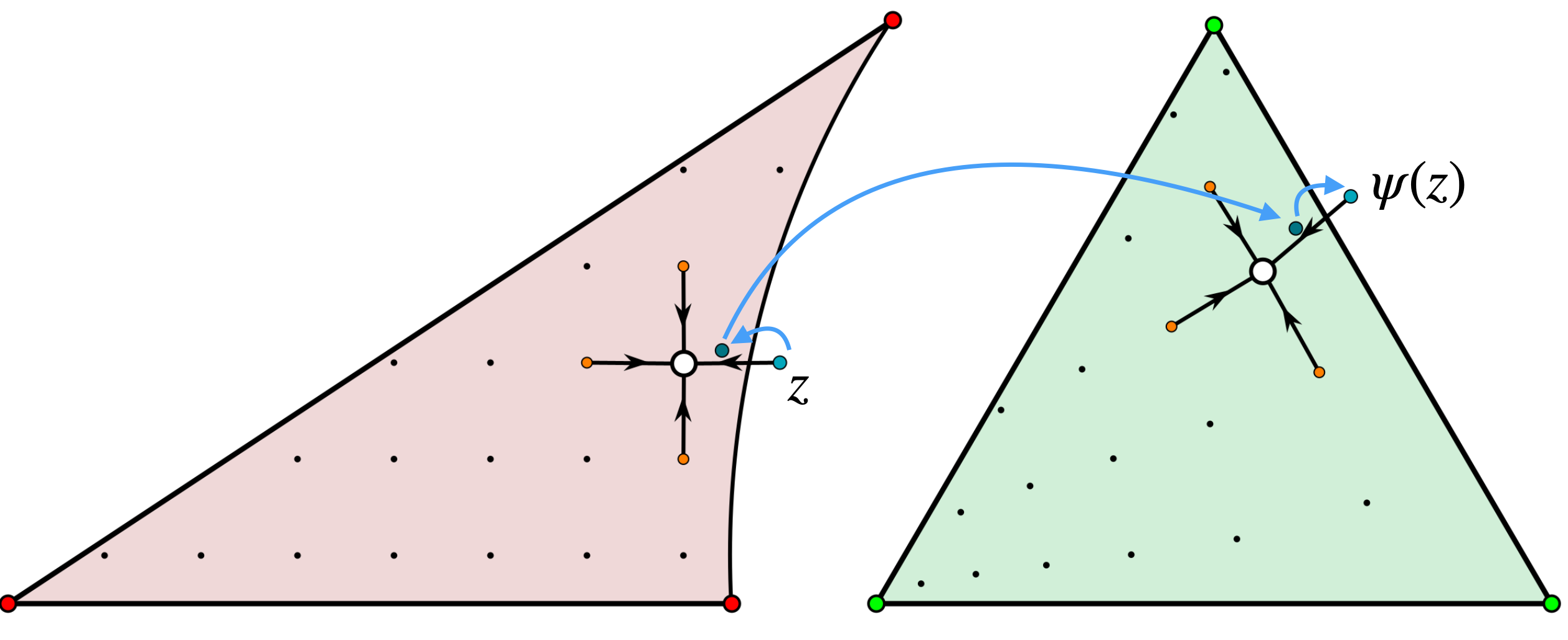}
		\caption{Applying the SRP on pixel grids}
	\label{applyingSRP}
%\end{figure}
\end{figure}

\begin{figure}[H]
\centering
\!\!\!\!\!\!\!\!\includegraphics[width=.52\textwidth]{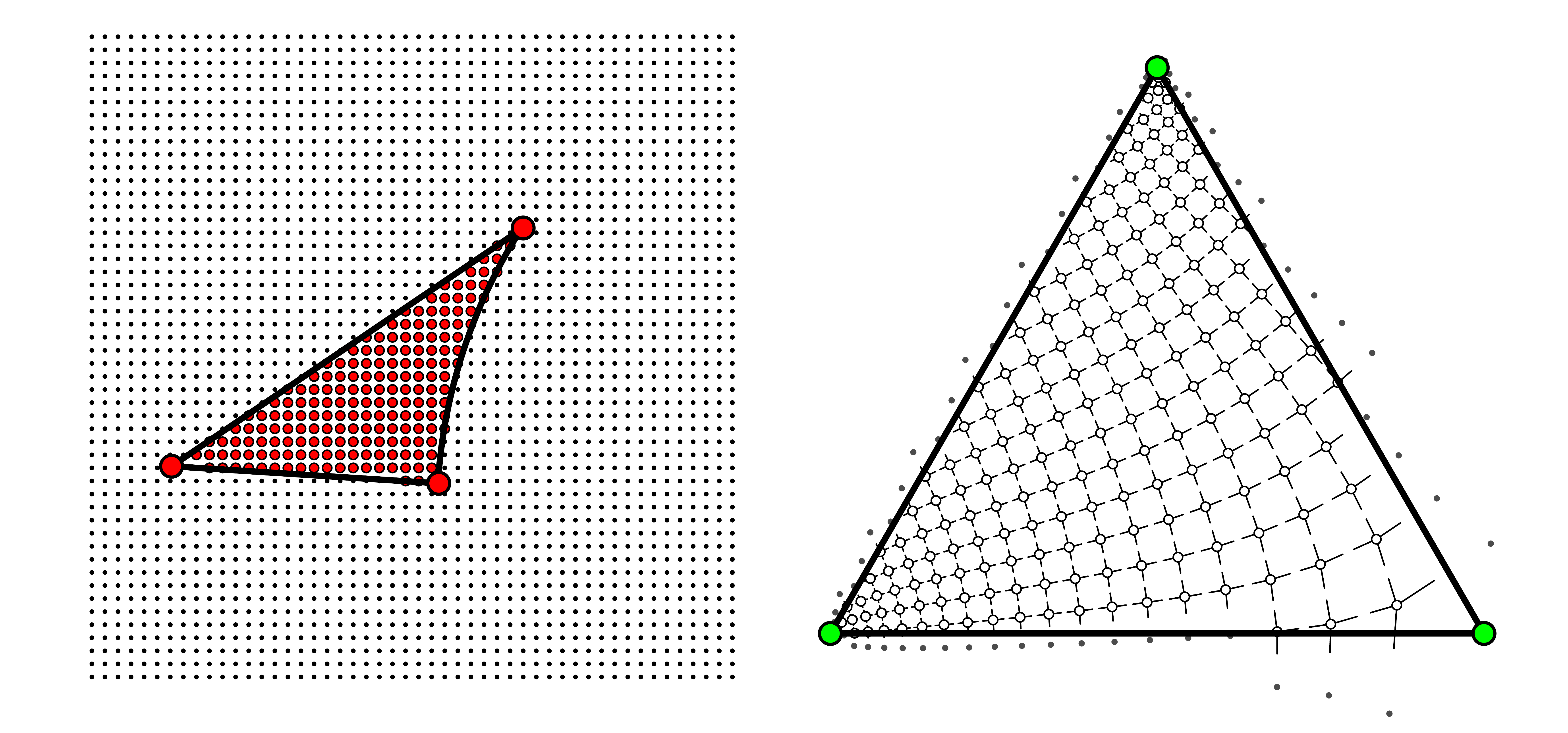}
\caption{The grid deformation}
\label{gridDeform}
\end{figure}

We explain the core idea of the algorithm by creating a conformal map $\psi \colon T_H \to T_E$ between hyperbolic and Euclidean triangles, respectively. As before, we choose this ``inverse'' direction as the colors get assigned to each pixel in the triangle $T_H$ by evaluating $f(\psi(z))$.
%The pixels form a square grid, and for each pixel we want to find the corresponding position in $T_E$, from which we take the color. This means for a pixel position $z$ we want to calculate $f(\psi^{-1}(z))$ for the map $\phi(z):=\psi^{-1}(z)$ described in the last section. The main problem is to find $\phi(z)$, which is conformal as well.
On a small scale, conformal maps behave like similarity transformations. So locally, small square grids are mapped to small square grids again. In a square grid, each point is in the center of gravity of its four neighbors. To approximate this behaviour in the limit case, we impose the same condition on the images of the pixels under~$\psi$.
We first describe the situation of the desired result, and then outline a parallelized algorithm to achieve it. Refer to Figures~\ref{applyingSRP}~and~\ref{gridDeform} for visual context of the following discussion.

We have to compute the map for every pixel position inside~$T_H$ as the resulting values describe the image. The pixels on the screen form a square grid with pixel distance~$\delta$.
For positions~$z$ whose four pixel neighbors lie inside~$T_H$, it is reasonable to assume that we will have
\[
\psi(z) = \frac{\psi(z+\delta) + \psi(z-\delta) + \psi(z+i\delta) + \psi(z-i\delta)}{4}
\]
in the discrete target situation.
In other words, the image of $z$ lies at the center of gravity of the images of its neighbors. Similar ideas have also been used in \cite{KD}~and~\cite{Polke}, though with weaker boundary conditions.

\columnbreak
At first sight, a problem seems to arise if not all of the pixel neighbors lie within~$T_H$, which typically happens very close to the boundary. In fact, this is exactly the place where we encode the geometry of the situation.
Suppose that $p$ is a pixel position near the boundary such that one of the four neighbors is not contained in $T_H$ as in Figure~\ref{applyingSRP}.
In this case, the average of the four neigboring images cannot be directly computed.
%as not all of them are precalculated.
To determine the missing images, we apply the SRP to find the analytic continuation.
For instance, if $z\not\in T_H$, we reflect it across the violated boundary edge to arrive at $z' \in T_H$. Now, its image $\psi(z') \in T_E$ can be found via linear interpolation of our other data points. Reflecting this positon with respect to the corresponding boundary edge in $T_E$ yields $\psi(z)$.

Finding the final configuration requires solving a large system of equations.
However, the final state can be interpreted as an energy minimum. Without going into detail, we can approximate this equilibrium by iteratively replacing each pixel’s position by the average of its neighbors. In the next section we will see that the problem can, in fact, be formulated as being linear.
Notably, the simplicity of the algorithm and the independence of the calculations for each pixel make it ideally suited for parallelization and GPU execution. In particular, standard GPU languages provide built-in support for linear interpolation. The following lines of pseudocode outline the core of the algorithm, very close to the literal code used in our implementations. Iterated applications of the algorithm create an approximation of the desired conformal map.

%\penalty20000
\pagebreak
\begin{lstlisting}[basicstyle={\ttfamily \footnotesize}]
map(p) := {
  q = p;
  if (inAB(p)) {q = reflAB(p)};
  if (inAC(p)) {q = reflAC(p)};
  if (inBC(p)) {q = reflBC(p)};
  image = interpolatedPixelMap(q);
  if (inAB(p)) {image = imgReflAB(image)};
  if (inAC(p)) {image = imgReflAC(image)};
  if (inBC(p)) {image = imgReflBC(image)};
  return image; };

forall (p in activePoints) {
  pixelMap[p] = (  map(p + (1,0))
           	 + map(p + (-1,0))
           	 + map(p + (0,1))
           	 + map(p + (0,-1)) ) / 4 };
\end{lstlisting}

Figure~\ref{gridDeform} shows a moderately sized pixel grid. The associated ``ghost points'', computed via the Schwarz reflection principle, are shown as small black dots outside the equilateral triangle.

\section{Linearity and Convergence}\label{convergence}
We examine the convergence behaviour of the algorithm next. Here, we only sketch the main ideas; for a detailed technical treatment, we refer to \cite{MRRG}.

Let us again assume that we are modelling the algorithm to calculate the map from a hyperbolic triangle $T_H$ to a corresponding Euclidean triangle $T_E$, both embedded in $\mathbb{C}$.
Pixel positions in the hyperbolic picture will be identified with points of the grid of Gaussian integers $\mathbb{G}=\mathbb{Z}+i\mathbb{Z}$. The hyperbolic plane will be represented by a Poincaré disk of sufficiently large radius $r >\!\!> 1$.

In principle, imagine that we select a suitable set of grid points in and around the hyperbolic triangle $T_H$, denoted by $\mathcal{A} \subset \mathbb{G}$, which serves as an index set for a vector $p \in \mathbb{C}^\mathcal{A}$. Each entry $p_z$ then encodes the (complex) image position of a grid point $z \in \mathcal{A}$.
Bilinear interpolation allows us to reconstruct an approximate function $\psi \colon T_H \to T_E$ in the following way:
For a non-grid point $z=a+ib$, set $\lfloor z \rfloor=\lfloor a \rfloor +i\lfloor b \rfloor\in \mathbb{G} $ and define the corners of a square cell that contains $z$ by 
\[ C(z):=\left( \lfloor z \rfloor,\lfloor z \rfloor+1,\lfloor z \rfloor +i, \lfloor z \rfloor+1+i \right).\]
Within each such cell, $z$ can be interpolated by bilinear combination from its corners. The factors are
\[\lambda(z) := \left( (1-\lambda)(1-\mu), \lambda(1-\mu), (1-\lambda)\mu, \lambda\mu\right),\] 
where $\lambda=a-\lfloor a\rfloor$ and $\mu=b-\lfloor b\rfloor$.
This allows us to write
\[ z=\sum_{i=1}^4 (\lambda(z))_i\cdot(C(z))_i. \]
On the image side, we use the approximation factors and set:
\begin{equation}
\label{eq:interpolatePsi}
\psi(z) := \sum_{i=1}^4 (\lambda(z))_i\cdot p_{C(z)_i}.
\end{equation}
Such bilinear interpolations between pixels are default primitives on GPUs that are easily accessible on the level of texture sampling.

At first glance, the involvement of circle inversions seems to suggest that the algorithm produces a non-linear system of equations.
However, the opposite is the case. The specific geometry of the situation makes it possible to formulate the algorithm using affine iteration steps of the form
\[ p\mapsto\mathbf{A}\cdot p +\mathbf{v}, \]
where the matrix $\mathbf{A}$ and the vector $\mathbf{v}$ encode the geometry of the hyperbolization process.

The points $\mathcal{A}$ that actively take part in our iteration are a subset of $\mathbb{G}$ sufficiently large that each point of $T_H$ is contained in a square cell for interpolation:
\[\mathcal{A} := \bigcup\limits_{z\in T_H}C(z).\]
We will consider $\mathcal{A}=\mathcal{I}\cupdot\mathcal{B}$ as a disjoint union of \emph{inner points} $\mathcal{I}=T_H\cap \mathbb{G}$ and \emph{boundary points} $\mathcal{B}$. For each point $z\in\mathcal{A}$ we will also need its pixel neighbors
\[N(z):=\left\{ z+1, z-1,z+i,z-i \right\}.\]

We formalize one iteration step in the following way: The vector $p$ in 
$\mathbf{A}\cdot p +\mathbf{v}$ will consist of $2\cdot |\mathcal{A}|$ 
complex entries where the first $m=|\mathcal{A}|$ will be indexed by the the points in $\mathcal{A}$. The following $m$ entries of $p$ are the complex conjugates of the first $m$ entries and will formally be indexed by a disjoint copy $\widetilde{\mathcal{A}}:=\{\widetilde{z} \mid z \in \mathcal{A} \}$ of $\mathcal{A}$.
We will organize our iteration step in a way such that the values of $p_z$ and $p_{\widetilde{z}}$ are complex conjugates. The reason for introducing the complex conjugates alongside the original coordinates is that the iteration involves Euclidean transformations that are either reflections or rotations. Per se, reflections are not linear operations over $\mathbb{C}$; they need complex conjugation. If we introduce an extra variable for each conjugated value, we can express reflections as linear operations –  now using the conjugates as input.

We will describe the update step for the points indexed by ${\mathcal{A}}$ next. The corresponding update steps for $\widetilde{\mathcal{A}}$ have to be defined similarly with interchanged roles of original and conjugate.
In each of the iteration steps, $p_z$ will be replaced by the average of its neighbors:
\[p_z \mapsto \frac{1}{4} \sum_{w\in N(z)} \mu(w),\]
where $\mu(w)$ corresponds to the {\tt map(...)} function of the above pseudocode. It reads out $p$ at the given position for inner points, and applies the SRP beforehand for points $w\not\in\mathcal{I}$:
\[
\mu(w)=
\begin{cases}
p_w&\mathrm{if }\ w\in \mathcal{I},\\
{\beta(w)}&\mathrm{otherwise},\\
\end{cases}
\]
where $\beta$ denotes the action of the SRP which operates as follows. Let $G_H$ and $G_E$ be the hyperbolic and Euclidean reflection groups, respectively, with fundamental cells $T_H$ and $T_E$. For any point $w\in\mathbb{G}$, there exists a unique symmetry operation $g\in G_H$ such that $g(w)\in T_H$.
To apply the SRP, we must find a corresponding symmetry $g' \in G_E$ mirroring the action of $g \in G_H$. We label the three reflections across the edges of $T_H$ as $\alpha, \beta, \gamma$ and their Euclidean counterparts as $\alpha', \beta', \gamma'$ – they generate $G_H$ and $G_E$, respectively.

Since we chose the Poincaré disk radius $r >\!\!> 1$ and since the relevant points $w$ for which we will have to determine the reflections are in a 1-neighborhood of the boundary points $\mathcal{B}$, they are all ``close'' to the boundary of $T_H$.
For most of these points, \emph{one} reflection suffices to move them into $T_H$. Points requiring more than one reflection will be near the corners of $T_H$ and only demand two of the reflections applied in alternating order a certain number of times.
For each such point $w$, we define $g_w$ as the reflection sequence of minimal length (involving at most two generators) such that $g_w(w) \in T_H$, an let $g'_w$ be the corresponding sequence in $G_E$.
The SRP is then modelled by
\begin{equation}
\beta(w):=(g'_w)^{-1}(\psi(g_w(w))).
\end{equation}
In other words, we first reflect $w$ via $g_w$ into $T_H$ and use equation \ref{eq:interpolatePsi} to bilinearly interpolate the image point $\psi(g_w(w))$ from the corners of the cell $C(g_w(w))$.
This is a linear combination of entries in~$p$! We then apply the inverse transformation $(g'_w)^{-1}$ on the Euclidean side to obtain $\beta(w)$.
The only potentially non-linear step is this final transformation which is either
\begin{itemize}
\item a Euclidean reflection $z\mapsto e^{it} \overline{z}+v$, or
\item a Euclidean rotation $z\mapsto e^{it} {z}+v$.
\end{itemize}

To accomodate both cases in a linear system, we use the complex conjugated points present in the vector~$p$. Rotations access the original values in $p$ (the first $m$ entries), while reflections access the conjugated part (the last $m$ entries).
In either case, the update equation for $p_w$ (and for $p_{\widetilde{w}}$) is an \emph{affine} function in the remaining entries of $p$, where the affine part comes from the shift term $v$.
Representing $p$ as $p^T=(q,\widetilde{q})^T$, we can write the iteration in block matrix form:
\begin{equation}
\label{eq:blockMatrix}
%\left( \begin{array}{c} \boxed{q}\\[1mm] \boxed{\widetilde{q}}\\ \end{array} \right) 
\begin{pmatrix} q \\ \widetilde{q} \end{pmatrix}
\mapsto
%\left( \begin{array}{cc} \boxed{M}&\boxed{N}\\[1mm] \boxed{\overline{N}}&\boxed{\overline{M}}\\ \end{array} \right)
\begin{pmatrix} M & N \\ \overline{N} & \overline{M} \end{pmatrix}
\cdot
%\left(\begin{array}{c} \boxed{q}\\[1mm] \boxed{\widetilde{q}} \\ \end{array} \right)
\begin{pmatrix} q \\ \widetilde{q} \end{pmatrix}
+
%\left( \begin{array}{c} \boxed{v}\\[1mm] \boxed{\overline{v}} \\ \end{array} \right)
\begin{pmatrix} v \\ \overline{v} \end{pmatrix}
\end{equation}
If the iteration is started with a vector $p$ whose last $m$ entries are the complex conjugates of the first $m$, then this property will be preserved throughout the iteration. The fixed point of the iteration process will therefore decompose into a component in the original space and a complex conjugate copy of it.

Each row of the matrix encodes the averaging of four points that are either original coordinates or calculated via the SRP by interpolating four cell coordinates, followed by an isometry. This implies that the row sum of the absolute values of the entries is $1$ in every row. This observation facilitates proving convergence of the iteration to a fixed point as it relates to one key property of the block matrix $\mathbf{A}$ from equation \ref{eq:blockMatrix}:
the spectral radius $\rho(\mathbf{A})$. This real number is the largest absolute value of an eigenvalue of $\mathbf{A})$. If $\rho(\mathbf{A}) < 1$, the iteration is contractive thus convergent.
In our setting, this convergence criterion can be certified by a powerful theorem due to Helmut Wielandt (1950), see \cite{Wielandt}.

\begin{figure*}[t]
	\centering
	\includegraphics[width=.98\textwidth]{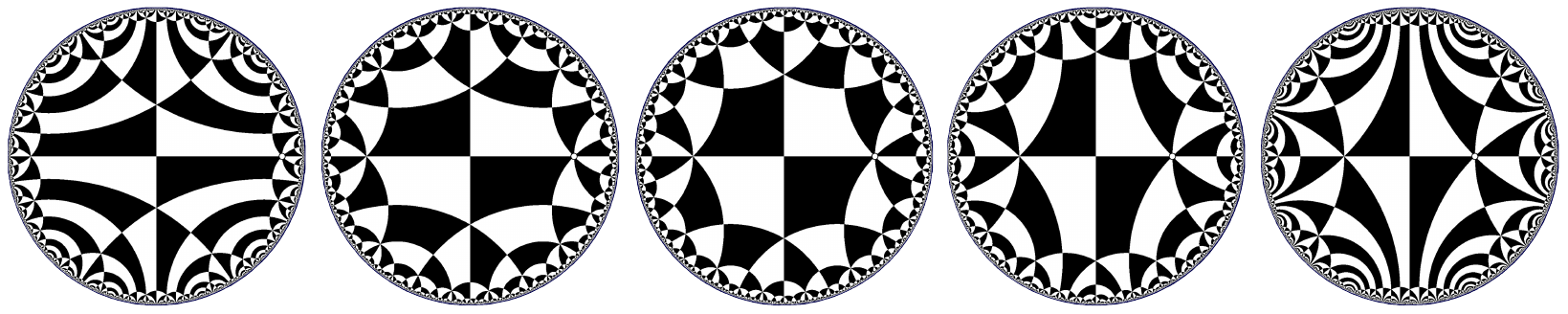}
	\caption{Same group, different conformal moduli}
	\label{modulus}
\end{figure*}
 
\begin{theorem}[Wielandt's Theorem]
Let $\mathbf{A} \in \mathbb{C}^{n\times n}$ and $\mathbf{B}\in \mathbb{R}_{\geq 0}^{n\times n}$ be two matrices.
If $\mathbf{B}$ is irreducible and $|\mathbf{A}_{i,j}| \leq \mathbf{B}_{i,j}$ for all $i,j \in [n]$, then $\rho(\mathbf{A}) \leq \rho(\mathbf{B})$.
Moreover, $\rho(\mathbf{A}) = \rho(\mathbf{B})$ holds if and only if there exist $\theta_0,\theta_1,\ldots,\theta_n \in \{ e^{it} \mid t \in \mathbb{R}\}$ with
\[ \mathbf{A}=\theta_0 \mathbf{D}^{-1}\mathbf{B}\mathbf{D} \]
for the diagonal matrix $\mathbf{D} = \mathrm{diag}(\theta_1,\ldots,\theta_n) \in \mathbb{C}^{n\times n}$.
\end{theorem}

We omit the technical details of the theorem's application in our setting and sketch the major steps instead:
\begin{itemize}
\item[1.] Restrict to the irreducible component of $\mathbf{A}$ that contains all inner points.
\item[2.] Apply Wielandt's theorem with $\mathbf{B}=|\mathbf{A}|$.
\item[3.] Note that $\rho(|\mathbf{A}|)\leq1$, since $||\mathbf{A}||_\infty=1$.
\item[4.] Show that there are no $\theta_i$ satisfying the equality criterion in the theorem.
\end{itemize}
In this manner, it follows that the spectral radius of $\mathbf{A}$ satisfies
$ \rho(\mathbf{A}) < \rho(\mathbf{|A|}) \leq 1 $
and thus convergence is established. The argument is fully explained in \cite{MRRG}.

\section{The Lower Symmetry Groups}
\label{moreover}

So far, we have studied the hyperbolization of Euclidean symmetry groups that are triangular reflection groups: $\mathbf{\ast 333}$, $\mathbf{\ast 442}$, and $\mathbf{\ast 632}$. Creating hyperbolizations of the remaining 14 wallpaper groups becomes increasingly difficult the less symmetries are present in the symmetry group. We briefly describe the situation in each case.

\subsection{\texorpdfstring{$\mathbf{333}$, $\mathbf{442}$, $\mathbf{632}$, $\mathbf{3{\ast}3}$, $\mathbf{4{\ast}2}$}{333, 442, 632, 3ast3, 4ast2}}
\label{sec:triangext}
 
These groups are closely related to the triangle reflection groups; they essentially live on the same geometric grids. More formally, if $G$ is a concrete embedding of one of these groups in the Euclidean plane, then one can add another generator to obtain a supergroup which is a triangle reflection group. In other words, one can draw a picture in such a group that \emph{looks like} belonging to a triangle reflection group. Each of these groups has index 2 in a corresponding $\mathbf{{\ast}ijk}$ group.

\begin{figure}[H]
	\centering
	\includegraphics[width=.22\textwidth]{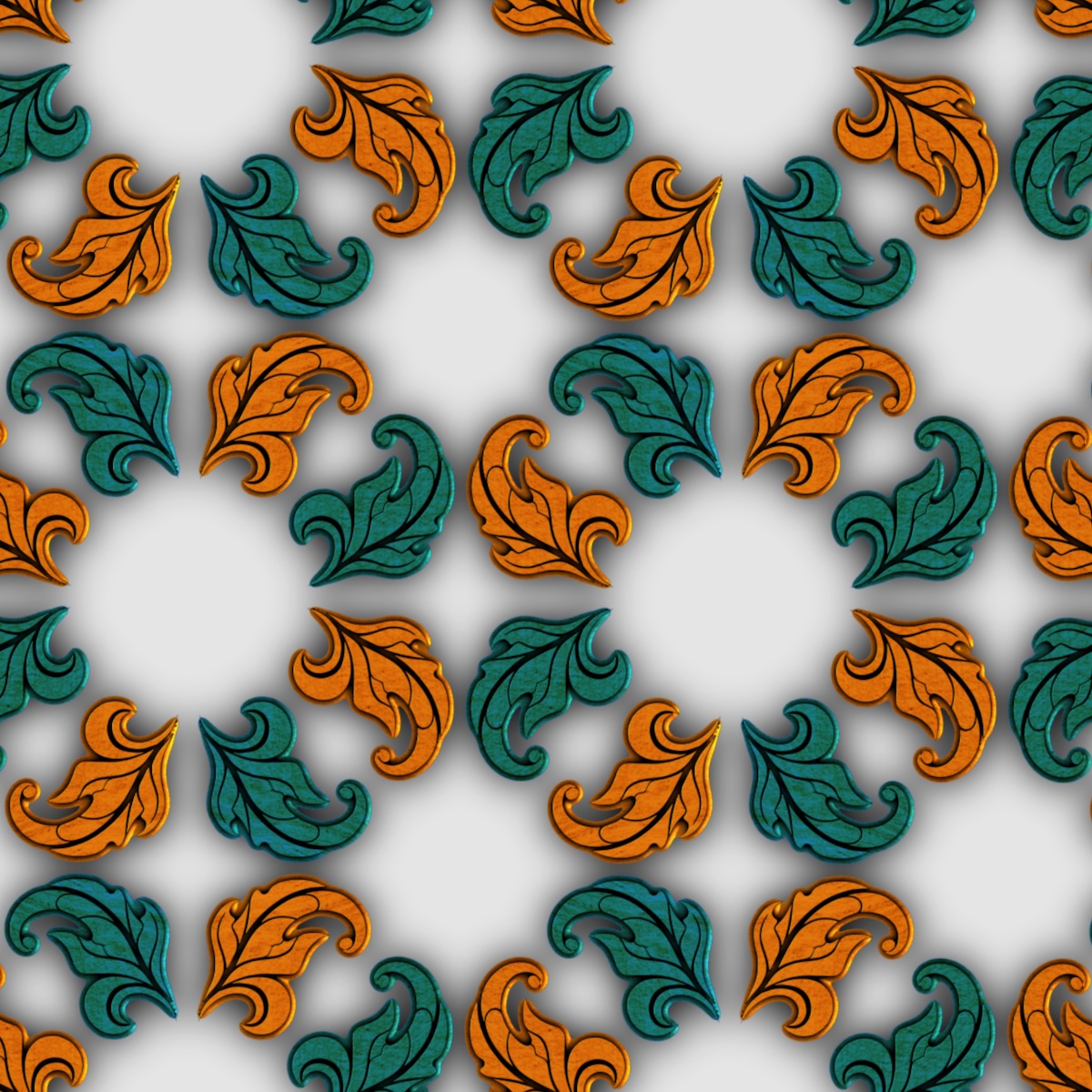}\quad
	\includegraphics[width=.22\textwidth]{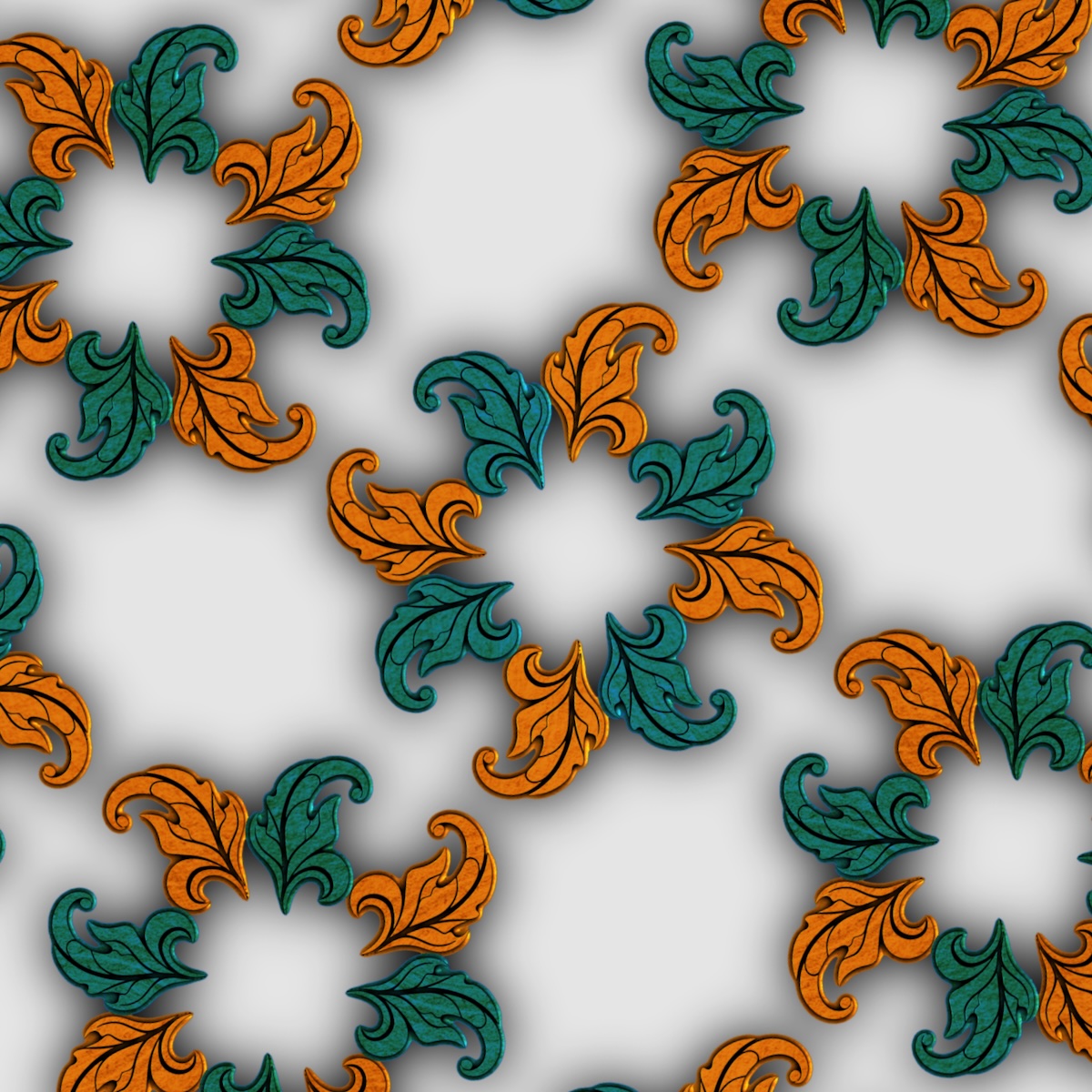}
	\caption{The relation of $\mathbf{442}$ and $\mathbf{4{\ast}2}$ to $\mathbf{{\ast}442}$}
	\label{highsym}
\end{figure}

One way to understand this relationship is via color symmetries. In Figure~\ref{highsym}, a pattern is shown in orange. In the left picture, the orange part is of type $\mathbf{442}$; in the right, it is of type $\mathbf{4{\ast}2}$. In both cases, a petrol part is added to the picture such that both parts together form a $\mathbf{{\ast}442}$ group.
This allows us to base a hyperbolization on the triangle group: it suffices to join two cells of the triangle group into a larger one that serves as fundamental region of the desired subgroup. Details are explained in \cite{HYP}. This method generally applies when the group in question is geometrically a subgroup of another group for which a hyperbolization has already been constructed.

\subsection{\texorpdfstring{$\mathbf{{\ast}2222}$}{ast2222}}

The next significant hurdle comes with $\mathbf{{\ast}2222}$, which is also a pure reflection group. In the Euclidean plane, it is generated by reflections across the sides of a rectangle $Q_E$.
At each corner of the rectangle, a 2-fold rotation center arises.
As before, increasing the index of any of the rotation centers yields a hyperbolization of the fundamental region. In this case, the resulting quadrilateral $Q_H$ may have different angels at each corner. When seeking a conformal transfer of the symmetry group, the fact that quadrilaterals have \emph{four} corners becomes relevant.
The Riemann mapping theorem guarantees the existence of conformal maps from both open quadrilaterals $Q_H$ and $Q_E$ to the unit disk. Each mapping is unique up to a Möbius automorphism of the disk and extends continuously to boundaries. 
Using this freedom, it is possible to map three of the corners of either quadritateral to the same positions.
However, the fourth vertex might or might not coincide. This depends on the \emph{conformal modulus} of the quadrilaterals: the cross ratio of the four image points on the circle.
This value is a conformal invariant of the quadrilateral and is a shape parameter similar to the aspect ratio of a rectangle. Figure~\ref{modulus} illustrates several hyperbolic reflection groups, each based on a fundamental quadrilateral with angles $(60^\circ,90^\circ,60^\circ,90^\circ)$ but with differing conformal moduli.
Thus, conformal hyperbolization of $\mathbf{{\ast}2222}$ depends on selecting the \emph{right} representation from a one-parameter family of combinatorially equivalent symmetry types. If this parameter is chosen correctly – which is, e.g., algorithmically possible via energy minimization – the Neandertal method works as before and produces a conformal map from $Q_H$ to $Q_E$. Extending this map via the SRP completes the hyperbolization for $\mathbf{{\ast}2222}$.

\subsection{\texorpdfstring{$\mathbf{{\ast}{\ast}}$, $\mathbf{{\ast}{\times}}$, $\mathbf{{\times{\times}}}$, $\mathbf{2{\ast}22}$, $\mathbf{22{\ast}}$, $\mathbf{22{\times}}$}{astast, astX, XX, 2ast22, 22ast, 22X}}
 
These are symmetry groups with even less symmetry sharing a common structural feature: In any geometric realization, each is a subgroup of a suitably chosen $\mathbf{{\ast}2222}$. In other words, in each of these groups a picture with rectangular reflection symmetry can be created. This allows to directly apply the same techniques as in Section~\ref{sec:triangext}.

\subsection{\texorpdfstring{$\mathbf{2222}$ and $\bigcirc$}{2222 and 0}}

The next problem arises for $\mathbf{2222}$, the group with four non-transitive rotation centers of order 2. Here, in general, two real parameters need to be adjusted for hyperbolization.
The fundamental region for this group can be any parallelogram. If the parallelogram is a rectangle, the methods from the previous section can again be applied. In the general case of a proper parallelogram, a second parameter has to be adapted to allow the creation of a conformal mapping.

\begin{wrapfigure}{R}{0.22\textwidth}
\centering
\includegraphics[width=0.22\textwidth]{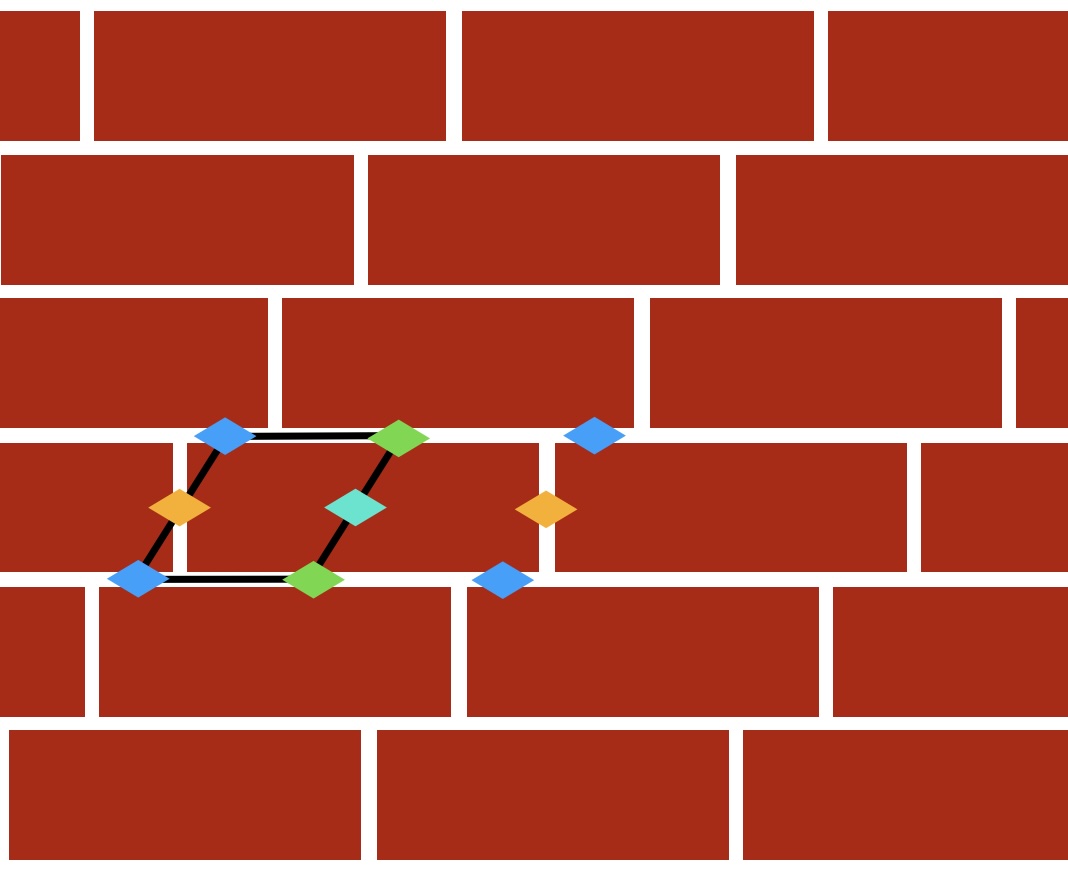}
\caption{A brick wall}
\label{fig:wall}
\end{wrapfigure}

One way to visualize this additional degree of freedom in the Euclidean case is via the analogy of a brick wall: While the bricks have an aspect ratio, every row of bricks additionally may be shifted relative to the previous one.
This shift is closely related to \emph{Dehn twists} in geometric group theory. In Figure~\ref{fig:wall}, a brick wall is shown with the corresponding fundamental region indicated as a black parallelogram.
Notice that the fundamental region only covers half of a brick, since one of the rotation centers (the cyan one)
lies inside the  brick.
All rotation centers of the ornament lie on the boundary of the parallelogram.

\begin{figure}[H]
	\centering
	\includegraphics[width=.47\textwidth]{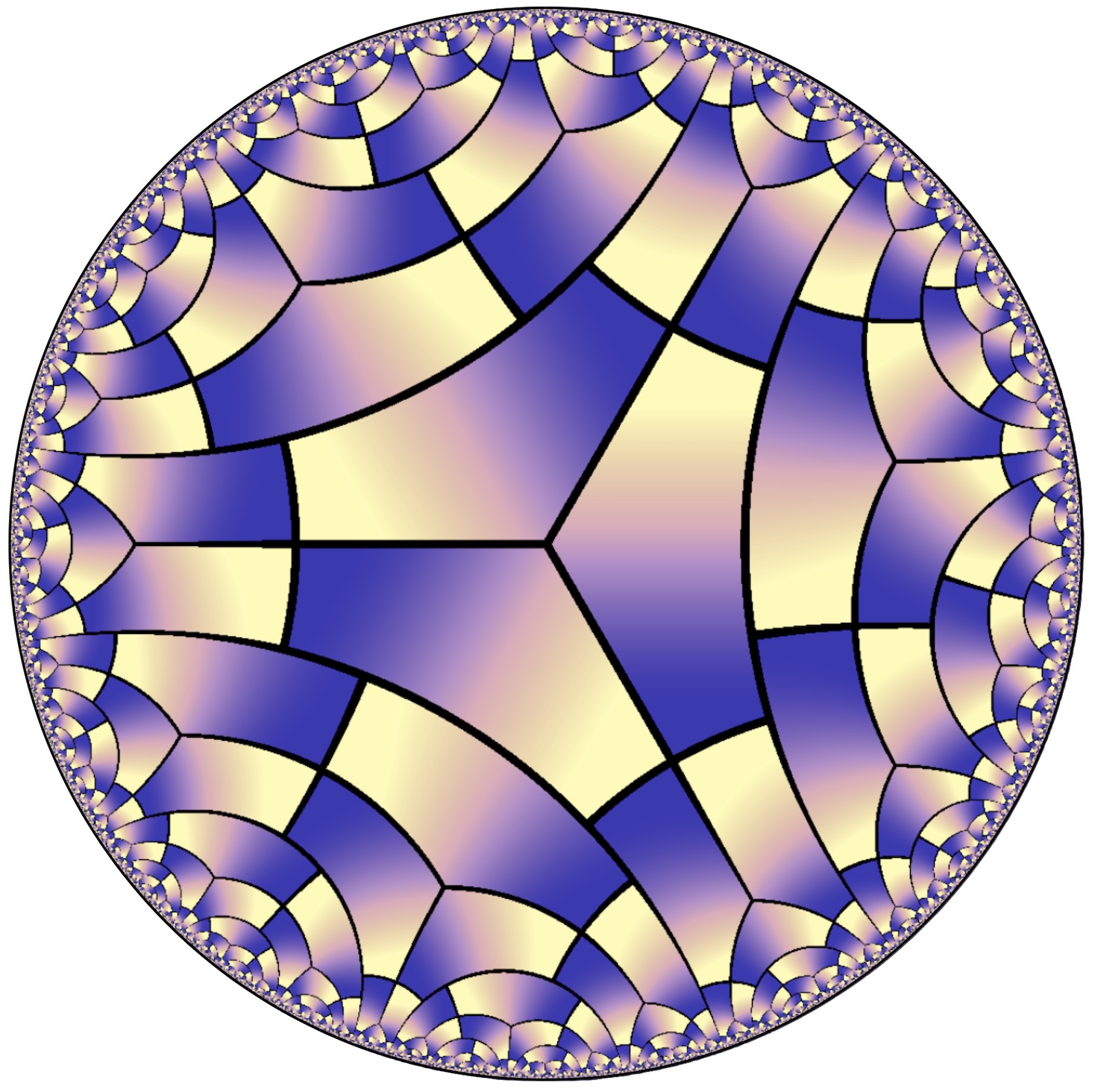}
	\caption{A hyperbolic brick wall}
	\label{fig:brick}
\end{figure}

Figure~\ref{fig:brick} shows the hyperbolic analogue of a brick wall. 
In this example the rotation center inside the Euclidean brick has been turned from 2-fold to 3-fold.
Thus the hyperbolic brick consists of 3 copies of the fundamental cell.
Each fundamental cell has 5 sides and the ``hyperbolic brick'' has 6 sides. Such a region can be seen in the center of
Figure~\ref{fig:brick}.
%\textcolor{red}{Here, each cell represents a fundamental region, i.e.\ a third of a brick.} In the given example, the fundamental region becomes a pentagon as one of the rotation centers along the edges has been modified to support a 3-fold rotation.

In this case, the Neandertal algorithm requires further adaption, as the underlying symmetry group lacks reflections and the SRP cannot be applied.
Nevertheless, it appears that a suitable generalization of the SRP can produce a unique conformal hyperbolization in this setting. Although no formal proof exists so far, numerical evidence looks promising. As before, if the two parameters are adjusted correctly, the Riemann mapping theorem guarantees the existence of a conformal map between the hyperbolic and the Euclidean fundamental cells.

The question remains how to extend this conformal map beyond the boundary of the fundamental region. Recall the discussion of the two approaches to such an extension in Section~\ref{sec:ConfAndSRP}: \emph{extension by symmetry} and \emph{extension by analytic continuation}.

For reflection groups, these two concepts coincided due to the SRP. We strongly conjecture that this is also the case for other symmetry groups. Our algorithm is exactly formulated to model \emph{extension by symmetry} for the boundary points.
This approach can be extended to other symmetry groups, enabling hyperbolizations of $\mathbf{2222}$ as well. If our conjecture holds, these hyperbolizations will automatically be conformal.
Finally, the group $\bigcirc$ – which contains only translations – can be treated as an index-2 subgroup of $\mathbf{2222}$ and thus hyperbolized accordingly.

\section{No Gritted Teeth}

%The remaining nine symmetry groups require more work but can be treated with similar methods. Each has at least one shape parameter and the hyperbolic groups have to be chosen consistently such that the conformal moduli of fundamental cells match.
%

Our methods are implemented in the ornament drawing program iOrnament~\cite{iOR}. They enable the user to switch between Euclidean and hyperbolic ornaments with ease. In particular, it is possible to hyperbolize own creations and to draw hyperbolic ornaments directly in the Poincaré disk. It is even possible to import external images like Escher's Euclidean tilings and create hyperbolizations of them.
The iterative nature of the algorithm has a beneficial side effect. As the map gradually approximates an equilibrium of forces, it can be used to animate smooth and visually pleasing transformations between symmetry groups.
A small online playground is available at \href{https://mathvisuals.org/Hyperbolization}{\tt mathvisuals.org/Hyperbolization}
~\cite{MV}.
\begin{figure}[H]
	\centering
	\includegraphics[width=.47\textwidth]{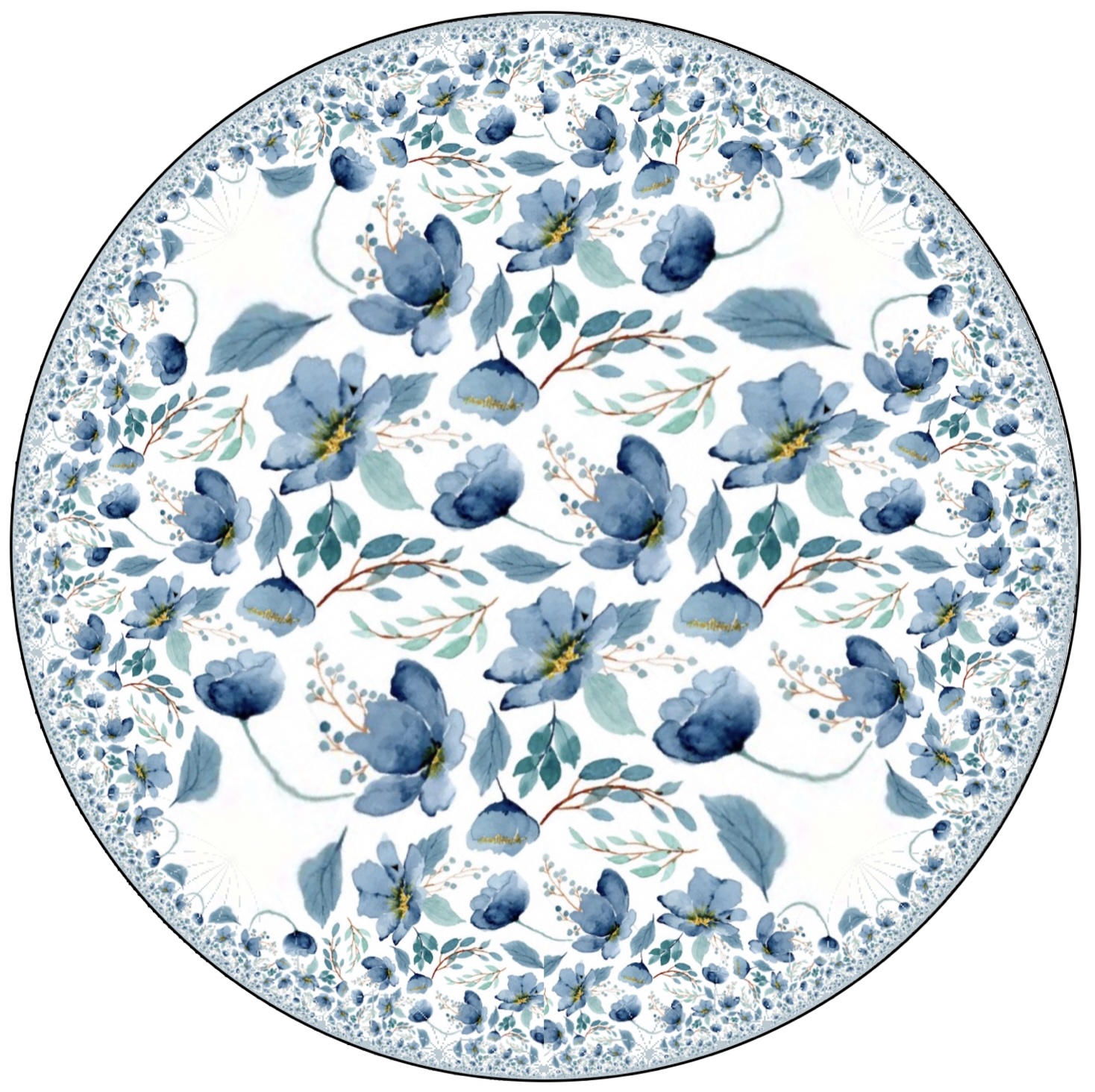}
	\caption{Hyperbolization of a $\bigcirc$ ornament}
	\label{fig:P1}
\end{figure}

One might wonder if the precision of the algorithm yields results that are sufficiently accurate. Measurements of local conformality in the fundamental cell indicate that we are extremely close to the desired results, optically indistinguishable.
Since the fundamental cell is usually large and since the iterated hyperbolic isometries are numerically stable, the pixel accuracy improves towards the boundary. By this, we can generate pixel maps up to 10000$\times$10000 resolution – the highest experimentally tested so far, with stunning sharpness even in large-scale prints.

\begin{figure*}[h]
\centering
\includegraphics[width=1\textwidth]{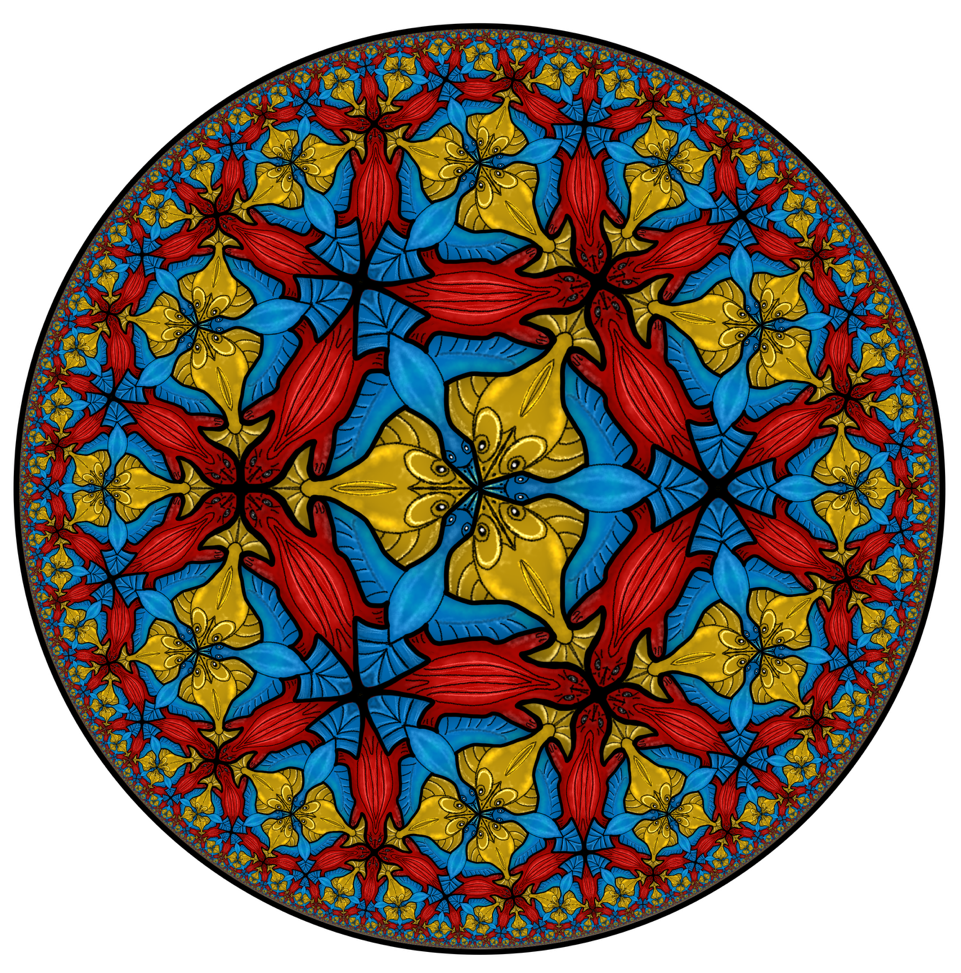}
\caption{A hyperbolic tiling, inspired by Escher's ``Fish, Duck, Turtle'' (drawn by the third author without gritted teeth)}
\vspace{10cm}
\label{EscherBFT}
\end{figure*}

We conclude our collection of examples with 
two hyperbolizations created directly in iOrnament. One of them Figure~\ref{fig:P1}
is a floral $\bigcirc$ ornament created directly from a wallpaper like pattern. The other one shown in 
 Figure~\ref{EscherBFT} is
our hand drawn take on the ``Fish, Duck, Turtle'' theme of M.C. Escher.

{\setlength{\baselineskip}{13pt} % tighten line spacing for bibliography
\raggedright % no right justification for References

} % end setlength, raggedright

\noindent All artworks by Escher are copyrighted by the M.\ C.\ Escher Company.
\end{multicols}

\end{document}